\documentclass[11pt]{article}
\usepackage{amsfonts,amssymb,amsmath}\usepackage{verbatim}
\def\sr{\stackrel}
\usepackage[dvips,matrix,curve,arrow]{xy}
\newtheorem{theorem}{Theorem}[section]
\newtheorem{lemma}[theorem]{Lemma}
\newtheorem{proposition}[theorem]{Proposition}
\newtheorem{corollary}[theorem]{Corollary}
\newtheorem{definition}[theorem]{Definition}

\newenvironment{proof}{\noindent\textsc{Proof:}\ }{\hfill$\square$\vskip1em}
\renewcommand{\geq}{\geqslant}
\newcommand{\Z}{{\mathbb Z}}\newcommand{\Q}{{\mathbb Q}}
\newcommand{\Qb}{{\overline{\mathbb Q}}}
\newcommand{\GL}{\operatorname{GL}}
\newcommand{\Hom}{\operatorname{Hom}}\newcommand{\End}{\operatorname{End}}
\newcommand{\Aut}{\operatorname{Aut}}\newcommand{\Ext}{\operatorname{Ext}}
\newcommand{\Inn}{\operatorname{Inn}}
\newcommand{\Br}{\operatorname{Br}}
\newcommand{\car}{\operatorname{char}}
\newcommand{\Gal}{\operatorname{Gal}}

\newcommand{\Res}{\operatorname{Res}}
\newcommand{\Inf}{\operatorname{Inf}}

\newcommand{\im}{\operatorname{im}}
\newcommand{\ab}{\operatorname{ab}}
\newcommand{\sign}{\operatorname{sgn}}
\newcommand{\Twist}{\operatorname{Twist}}

\title{Curves of genus 2 with group of automorphisms\\ isomorphic to $D_8$ or $D_{12}$}
\author{Gabriel Cardona$^*$
  \and
  Jordi Quer\thanks{Supported by BFM2000-0794-C02-02 and HPRN-CT-2000-00114}} 
\begin{document}\maketitle

\begin{abstract}
The classification of curves of genus 2 was studied by Clebsch and Bolza
using invariants of binary sextic forms,
and completed by Igusa with the computation of the three-dimensional
moduli variety $\mathcal M_2$ classifying genus 2 curves over an algebraically closed field.
The locus of curves with group of automorphisms isomorphic
to one of the dihedral groups $D_8$ or $D_{12}$ is a one-dimensional subvariety.

In this paper we classify these curves over an arbitrary perfect field $k$
of characteristic $\car k\neq2$ in the $D_8$ case and $\car k\neq2,3$ in the $D_{12}$ case.
We first parametrize the $\overline k$-isomorphism classes of curves defined over $k$
by the $k$-rational points of a quasi-affine one-dimensional subvariety of $\mathcal M_2$;
then, for every curve $C/k$ representing a point in that variety we compute all its $k$-twists,
which is equivalent to the computation of the cohomology set $H^1(G_k,\Aut(C))$.

The classification is always performed by explicitly describing the objects involved:
the curves are given by hyperelliptic models
and their groups of automorphisms represented as subgroups of $\GL_2(\overline k)$.
In particular, we give two generic hyperelliptic equations,
depending on several parameters of $k$,
that by specialization produce all curves in every $k$-isomorphism class.

As an application of the classification of curves of genus 2 obtained,
we get precise arithmetic information on their elliptic quotients and on their jacobians.
Over the field $k=\Q$, we show that the elliptic quotients of the curves with automorphisms
$D_8$ and $D_{12}$ are precisely the $\Q$-curves of degrees 2 and 3, respectively,
and we determine which curves have jacobians of $\GL_2$-type.
\end{abstract}

\section{Preliminars on hyperelliptic curves and curves of genus 2}

This section contains basic definitions, notation,
and some well-known facts on hyperelliptic curves and curves of genus 2.
References are \cite{prolegomena}, \cite{mestre}, \cite{poonen}.

Thorough the paper $k$ is a perfect field of characteristic different from 2,
and $G_k$ is the Galois group of an algebraic closure $\overline k/k$.
The Galois action on the elements of any $G_k$-set will be denoted exponentially on the left:
$(\sigma,a)\mapsto{}^\sigma a$ for $\sigma\in G_k$ and $a$ in a $G_k$-set.

Some results are stated in terms of elements of $\Br_2(k)\simeq H^1(G_k,\{\pm1\})$,
the 2-torsion of the Brauer group of the field $k$.
We denote by $(a,b)$ the class of the quaternion algebra with basis $1,i,j,ij$
and multiplication defined by $i^2=a,j^2=b,ji=-ij$.

A curve $C/k$ of genus $g$ is hyperelliptic over $k$
if there is a morphism $C\to\mathbb P^1$ of degree 2 defined over $k$.
These curves have a model given by an equation
\begin{equation}\label{hyperelliptic_equation}
C:\quad Y^2=F(X)
\end{equation}
with $F(X)\in k[X]$ a polynomial of degree $2g+1$ or $2g+2$ without multiple roots.
Conversely, every such equation is a model of a curve of genus $g$
defined over $k$ and hyperelliptic over $k$,
with a unique singularity at the point at infinity
that corresponds to one or two points in a regular model,
depending on whether the degree of the polynomial $F(X)$ is odd or even.
We call the equations (\ref{hyperelliptic_equation}) \emph{hyperelliptic equations}.

Every curve $C/k$ of genus 2 is hyperelliptic over $k$.
In general, we will implicitly assume that curves of genus 2 over $k$
are given by hyperelliptic equations.
The isomorphisms between two such curves correspond,
in terms of hyperelliptic equations, to transformations of the type
\begin{equation}\label{isomorfismes}
X'=\frac{aX+b}{cX+d},\qquad Y'=\frac{(ad-bc)Y}{(cX+d)^3},\qquad
   M=\begin{pmatrix} a & b \\ c & d \end{pmatrix}\in\GL_2(\overline k),\end{equation}
and the field of definition of that isomorphism is the field generated
over the common field of definition of the two curves by the coefficients of the matrix $M$.
In particular, for every curve $C/k$ of genus 2, the group of automorphisms $A=\Aut(C)$
can be identified with a subgroup of $\GL_2(\overline k)$ which is closed by
the Galois action of the group $G_k$.

\paragraph{Twisting curves.}
For every curve $C$ defined over a field $k$,
the set of its twists $\Twist(C/k)$
is the set of $k$-isomorphism classes of curves $C'/k$
that are $\overline k$-isomorphic to $C$.
The map sending an isomorphism $\phi:C'\to C$ to the class of the 1-cocycle
${}^\sigma\phi\circ\phi^{-1}:G_k\to\Aut(C)$
is a bijection between the set of $k$-twists of $C$
and the cohomology set $H^1(G_k,\Aut(C))$.

For curves $C/k$ of genus 2 given by hyperelliptic equations
the inverse map is easily computed from the identification
of $\Aut(C)$ with a subgroup of $\GL_2(\overline k)$.
Indeed, given a 1-cocycle $\sigma\mapsto\xi_\sigma$ of $G_k$ with values in $\Aut(C)$
we may view it as taking values in $\GL_2(\overline k)$.
Since the cohomology set $H^1(G_k,\GL_2(\overline k))$ is trivial,
there is a matrix $M\in\GL_2(\overline k)$
such that $\xi_\sigma={}^\sigma M\cdot M^{-1}$.
Then the curve $C'$ obtained from $C$ and the matrix $M$ using (\ref{isomorfismes})
is the twist of $C$ corresponding to the cohomology class of the 1-cocycle $\xi$.

Every hyperelliptic curve $C/k$ has a canonical involution,
the \emph{hyperelliptic involution} $\iota$,
that is always defined over the ground field
and commutes with every automorphism.
In the matricial representation of $\Aut(C)$ for the genus 2
case it corresponds to the matrix $-1$
and is the unique nontrivial automorphism given by a scalar matrix.
The quotient group $A'=A/\langle\iota\rangle$ is known as the \emph{reduced group of automorphisms}.
The sequence $1\rightarrow \langle\iota\rangle\rightarrow A\rightarrow A'\rightarrow 1$
is an exact sequence of $G_k$-groups.
Since the group on the left is contained in the center of the group on the middle
one obtains (cf. \cite[Section 8.3]{serre})
 the following long exact sequence of cohomology sets
$$\dots\longrightarrow {A'}^{G_k}\buildrel\delta\over\longrightarrow H^1(G_k,\langle\iota\rangle)
  \longrightarrow H^1(G_k,A)\longrightarrow H^1(G_k,A')\buildrel\Delta\over\longrightarrow
  H^2(G_k,\langle\iota\rangle)$$
from which one deduces the exact sequence
\begin{equation}\label{cohomology_exact_sequence}
1\longrightarrow H^1(G_k,\langle\iota\rangle)/\delta({A'}^{G_k})\longrightarrow
  H^1(G_k,A)\longrightarrow H^1(G_k,A')[\Delta]\longrightarrow1.
\end{equation}
The group $H^1(G_k,\langle\iota\rangle)/\delta({A'}^{G_k})$
acts faithfully on the set $H^1(G_k,A)$ and the orbits of that action
are in bijection with the cohomology set $H^1(G_k,A')[\Delta]$ (cf. \cite[Section 8.4]{serre}).
Hence, there is a bijection between the set of twists of the curve $C/k$
and the set of pairs consisting of an element of
$H^1(G_k,\langle\iota\rangle)/\delta({A'}^{G_k})$
and an element of $H^1(G_k,A')[\Delta]$.

The group $H^1(G_k,\langle\iota\rangle)/\delta({A'}^{G_k})$ can be identified with a quotient of
the group $k^*/k^{*2}$ by a finite subgroup,
and the 1-cocycles representing its elements are homomorphisms $\xi:G_k\to\langle\iota\rangle$.
If $k(\sqrt d)$ is the fixed field of the kernel of this homomorphism,
then the curve obtained by twisting a curve with equation $Y^2=F(X)$
by the element of $H^1(G_k,A)$ corresponding to $\xi$
is the curve with equation $Y^2=dF(X)$ or, up to $k$-isomorphism, with equation $dY^2=F(X)$.
We call these twists \emph{hyperelliptic twists},
and we say that two twists of a curve differ by an hyperelliptic twist
if they correspond to two elements of $H^1(G_k,A)$ with the same image in the set $H^1(G_k,A')$,
which is equivalent to the fact that the two curves admit hyperelliptic models
(\ref{hyperelliptic_equation})
with polynomials $F(X)$ that differ only by the multiplication by some nonzero element of $k$.

It is immediate to check that if two curves differ by an hyperelliptic twist
then their groups of automorphisms are isomorphic as $G_k$-groups.
The converse is not in general true, but it is ``almost'' true
for the two families of hyperelliptic curves that we will study in this paper,
namely curves of genus 2 with groups of automorphisms
isomorphic to one of the dihedral groups $D_8$ or $D_{12}$.
Indeed, as we will see in section 4, for curves of genus 2 with $\Aut(C)\simeq D_{12}$
the Galois action on $\Aut(C)$ is enough for classifying twists up to hyperelliptic twist.
The situation in the case $\Aut(C)\simeq D_8$ is a bit more complicated:
a curve can have either only one or at most two twists, up to hyperelliptic twists,
with the same Galois action on $\Aut(C)$;
equivalently, the set of elements of $H^1(G_k,A)$ that produce twisted curves with
group of automorphisms isomorphic to a given $G_k$-group
can have either one or two images in the set $H^1(G_k,A')[\Delta]$.

\section{Parametrization over an algebraically closed field}

The complete classification of curves of genus 2 over an algebraically closed field
was obtained by Igusa in \cite{igusa}, completing previous work by Bolza and Clebsch.
The classification is obtained in terms of invariants of binary sextics
$$a_0Y^6+a_1Y^5X+a_2Y^4X^2+a_3Y^3X^3+a_4Y^2X^4+a_5YX^5+a_6X^6$$
attached to the polynomials $F(X)=\sum_{i=0}^6 a_iX^i$
that correspond to hyperelliptic equations of the curves.
An \emph{invariant of binary sextic forms} of degree $d\geq 1$
is a polynomial expression $I\in k[a_0,\dots,a_6]$ in the coefficients of the sextic
that after a linear transformation of the variables
changes by the $d$-th power of the determinant.

The invariants that normally are used in the litterature are the Clebsch invariants,
denoted as $A,B,C,D$ in \cite{bolza} and \cite{mestre},
that we will denote by $c_2,c_4,c_6,c_{10}$,
and the Igusa invariants $I_2,I_4,I_6,I_{10}$,
defined as symmetric expressions of the roots of $F$ (cf. \cite[pag. 620]{igusa}).
In both cases of our notation the subindices give the degrees of the invariants.
These invariants are not defined in charactersitic $2$
and Clebsch invariants do not behave nicely in characteristics $3$ and $5$.
Another set of invariants reducing well in all characteristics
was introduced by Igusa in \cite{igusa} (he called them ``arithmetic invariants''),
but for the purposes of this paper the $c_j$ and $I_j$ will be enough.

The \emph{absolute invariants} are defined as the quotients of invariants of the same degree.
The classification of curves of genus 2 up to $\overline k$-isomorphism is given by their
absolute invariants: two curves are isomorphic if, and only if,
they have the same absolute invariants.

The possible reduced groups of automorphisms of curves of genus 2 were determined by Bolza
in terms of their invariants (cf. \cite[pag. 70]{bolza}),
and the structure of the corresponding groups can be found in \cite{esquela}.
The picture, outside from characteristics 2, 3 and 5, is the following:
the group $\Aut(C)$ is isomorphic to one of the groups
$$C_2,\ \ V_4,\ \ D_8,\ \ D_{12},\ \ 2D_{12},\ \ \widetilde S_4,\ \ C_{10}$$
with $2D_{12}$ and $\widetilde{\mathfrak S}_4$ denoting certain double covers
of the dihedral group $D_{12}$ and the symmetric group $\mathfrak S_4$.
The moduli space of curves of genus 2 is of dimension 3.
The generic curve has group of automorphisms isomorphic to $C_2$.
The curves with $V_4\subseteq\Aut(C)$ cut out a surface in that 3-dimensional moduli space,
and those with $D_8\subseteq\Aut(C)$ or with $D_{12}\subseteq\Aut(C)$
describe two curves contained in that surface, intersecting in a unique point.
Each of the three remaining groups is the group of automorphisms of a unique curve up to isomorphism,
and the three corresponding moduli points are, respectively,
the intersection of the $D_8$-curve and the $D_{12}$-curve,
a point contained in the $D_8$-curve but not in the $D_{12}$ one,
and a point not contained in the $V_4$-surface.
Small modifications of the picture are needed in characteristics 3 and 5.

In this paper we study the two families of curves of genus 2
with $\Aut(C)\simeq D_8$ or $D_{12}$.
As explained in the previous paragraph,
over an algebraically closed field these curves are one-parametric families,
and it is not difficult to find explicit parametrizations (cf. \cite[pag. 50]{bolza}).
The next two propositions introduce parametrizations in terms of an absolute invariant.
Our parametrizations are different from those of \cite{bolza},
which are the usual in the litterature,
and have the advantage that the field of definition of the curve
is also the field of definition of the point of moduli.

\begin{proposition}\label{proposition_absolute_invariant_D8}
Let $k$ be a field of characteristic different from 2.
The $\overline k$-isomorphism classes of curves $C/k$ of genus 2
with $\Aut(C)\simeq D_8$ are clasified by the open subset of the affine line
$$k\smallsetminus\{0,1/4\}.$$
An explicit bijection is obtained by associating to a given curve $C$ the element
$$t=\begin{cases} \displaystyle\frac{I_4}{I_2^2},&\quad\text{if}\quad\car k=5,\\
&\\
\displaystyle\frac{2I_2^3}{I_6-I_2^3},&\quad\text{if}\quad\car k=3,\\
&\\
\displaystyle\frac{8c_6(6c_4-c_2^2)+9c_{10}}{900c_{10}},&\quad\text{otherwise}.\end{cases}$$
and to an element $t\in k\smallsetminus\{0,1/4\}$ the curve with equation
$$Y^2=X^5+X^3+t\,X.$$
\end{proposition}

\begin{proof}
The polynomial $X^5+X^3+t\,X$ has multiple roots for exactly the two values $t=0,1/4$,
hence the equation $Y^2=X^5+X^3+t\,X$ defines an hyperelliptic curve $C_t$ of genus 2
for every $t\neq0,1/4$.
One checks that the group of automorphisms of this curve is the group
isomorphic to $D_8$ generated by the two matrices
$$U=\begin{pmatrix} -\sqrt{-1} & 0 \\ 0 & \sqrt{-1} \end{pmatrix},\qquad
  V=\begin{pmatrix} 0 & 1 \\ 1 & 0 \end{pmatrix}.$$
Conversely, it is easily seen that every curve $C/\overline k$ with $\Aut(C)\simeq D_8$
is $\overline k$-isomorphic to a curve $C_t$ for some $t\in\overline k$.

Using the formulas expressing the Clebsch and Igusa invariants of a curve
as polynomials on the coefficients of a polynomial of degree 5 or 6
in an hyperelliptic equation defining the curve, one checks that $t\in\overline k$
is an absolute invariant of the curve $C_t$,
written in terms of Clebsch or Igusa invariants by the formulas given in the statement
of the proposition.

Since $t$ is an absolute invariant of the curve $C_t$ it follows that $C_t$ and $C_{t'}$
are isomorphic if, and only if, $t=t'$.

As for the field of definition, if $t\in k$ the curve $C_t$ is defined over $k$.
Conversely, if a curve $C$ is defined over $k$ and $C_t$ is a curve isomorphic to it,
then $t\in k$ since $t$ is an absolute invariant of $C$.
\end{proof}

The next proposition is analogous for $\Aut(C)\simeq D_{12}$ and the proof is the same.

\begin{proposition}\label{proposition_absolute_invariant_D12}
Let $k$ be a field of characteristic different from 2 and 3.
The $\overline k$-isomorphism classes of curves $C/k$ of genus 2
with $\Aut(C)\simeq D_{12}$ are clasified by the open subset of the affine line
$$k\smallsetminus\{0,1/4\}.$$
An explicit bijection is obtained by associating to a given curve $C$ the element
$$t=\begin{cases} \displaystyle\frac{-I_4}{I_2^2},&\quad\text{if}\quad\car k=5,\\
&\\
\displaystyle\frac{3c_4c_6-c_{10}}{50c_{10}},&\quad\text{otherwise}.\end{cases}$$
and to an element $t\in k\smallsetminus\{0,1/4\}$ the curve with equation
$$Y^2=X^6+X^3+t.$$
\end{proposition}

From now on we will use the invariants in these propositions
to parametrize the curves up to $\overline k$-isomorphism,
and for that we make the following

\begin{definition}
Let $C/k$ be a curve of genus 2 with $\Aut(C)\simeq D_8$ or $D_{12}$.
The element $t=t(C)\in k\smallsetminus\{0,1/4\}$ defined as a quotient of invariants
of the same degree by the formulas of the previous propositions
will be called \emph{the absolute invariant} of the curve $C$.
\end{definition}

\section{Galois structures on $D_8$ and $D_{12}$}

Let $A$ be a finite group.
The $G_k$-group structures on the group $A$
(up to isomorphism of $G_k$-groups) are classified by the cohomology set
$$H^1(G_k,\Aut(A))=\Inn(\Aut(A))\backslash\Hom(G_k,\Aut(A))$$
with $\Aut(A)$ viewed as a $G_k$-module with trivial action.

A group action $(\sigma,a)\mapsto{}^\sigma a:G_k\times A\to A$
corresponds to a homomorphism $\rho:G_k\to\Aut(A)$,
and isomorphic actions differ by conjugation by an element of $\Aut(A)$.
The fixed field of $\ker\rho$ will be called
the \emph{field of definition} of the $G_k$-action on $A$ (or of $A$ viewed as a $G_k$-group),
and from now on it will be denoted by $K$.

Let $A$ be a $G_k$-group.
From every 1-cocycle $\sigma\mapsto\xi_\sigma:G_k\to A$ we can define another
action on the group $A$ by the formula ${}_\xi^\sigma a=\xi_\sigma^{-1}{}^\sigma a\xi_\sigma$.
The corresponding $G_k$-group ${}_\xi A$ depends only on the cohomology class of $\xi$
and is known as the \emph{twisted group}.
Two curves defined over $k$ and isomorphic over $\overline k$
have Galois action on their groups of automorphisms that are twisted of each other.
A group action on $A$ induces an action on the maximal abelian quotient $A^{\ab}$,
and it is clear that twisted actions induce the same action on that quotient.

Given a conjugation class $T$ of subgroups of $\Aut(A)$,
we say that a $G_k$-group structure on $A$ is \emph{of type $T$} if $\im\rho\in T$.
Then the action induces an isomorphism $\Gal(K/k)\simeq\im\rho$
and two $G_k$-group structures belonging to the same type $T$
and having the same field $K$ as the field of definition,
with corresponding maps $\rho,\rho':G_k\to\Aut(A)$,
differ in an isomorphism $\im\rho\simeq\im\rho'$
up to isomorphisms induced by conjugation by an element of $\Aut(A)$.

\paragraph{Galois structures on $D_8$.}
Consider the group $A\simeq D_8$, with presentation
$$A=\langle\ U,V\ |\ U^2=V^4=1,VU=UV^3\ \rangle.$$
The group $\Aut(A)$ is also isomorphic to $D_8$ and it is generated by the two automorphisms
$s$ and $t$ defined by
\begin{equation}\label{automorphisms_D8}
\begin{aligned} s(U) &= U  \\ s(V) &= V^3 \end{aligned}\qquad
\begin{aligned} t(U) &= UV \\ t(V) &= V \end{aligned}
\end{equation}
with the relations $s^2=t^4=1,ts=st^3$.
The inner automorphisms of $\Aut(A)$ are the subgroup of $\Aut(\Aut(A))$,
isomorphic to the Klein group $V_4$,
generated by the conjugation by $s$ and the conjugation by $t$.

The characteristic subgroups of $A$ (the subgroups invariant by every automorphism) are,
apart from the trivial subgroup and the group $A$ itself,
the subgroups $Z(A)=\langle V^2\rangle\simeq C_2$ and $\langle V\rangle\simeq C_4$.
In particular, every $G_k$-group structure on $A$
induces corresponding structures on these groups
and also on the quotients $A'=A/Z(A)\simeq V_4$ and $A/\langle V\rangle\simeq C_2$.
We remark that $A'=A^{\ab}$.

For every $G_k$-group structure on $A$ with field of definition $K$
we will denote by $K_1$ the field of definition of the induced action on $A'$
and by $K_2$ the field of defintion of the induced action on the subgroup $\langle V\rangle$.
Then $\Gal(K/k)$ is isomorphic to a subgroup of $D_8$
and $\Gal(K_i/k)$ can be trivial or isomorphic to the group $C_2$ for $i=1,2$.
Through the rest of the paper,
whenever a $G_k$-group $A$ isomorphic to $D_8$ is considered,
$K,K_1$ and $K_2$ will always denote the extensions of $k$
that are the fields of definition of the Galois action on the group $A$,
on the quotient $A'$ and on the subgroup $\langle V\rangle$, respectively.

By examination of all the subgroups of $\Aut(A)$
one finds all the possible types of $G_k$-group structures,
that are classified in the following table

$$\begin{array}{|c|c|c|c|c|c|}
\hline
\im\rho & \Gal(K/k) & \Gal(K_2/k) & \Gal(K_1/k) & \text{Type} & \text{remarks} \\
\hline\hline
1 & I & I & I & I & \\
\hline
\langle t^2\rangle    & C_2 &   I &   I & C_2^A & \\
\hline
\langle s\rangle,\langle st^2\rangle  & C_2 & C_2 &   I & C_2^B & \\
\hline
\langle st\rangle,\langle st^3\rangle & C_2 & C_2 & C_2 & C_2^C & K_1=K_2 \\
\hline
\langle t\rangle      & C_4 &   I & C_2 & C_4 & (u,-1)=1 \\
\hline
\langle s,t^2\rangle  & V_4 & C_2 &   I & V_4^A & \\
\hline
\langle st,t^2\rangle & V_4 & C_2 & C_2 & V_4^B & K_1=K_2 \\
\hline
\langle s,t\rangle    & D_8 & C_2 & C_2 & D_8 & (u,-v)=1 \\
\hline
\end{array}$$
\begin{center}{\bf Table 1.} $G_k$-group structures on $D_8$\end{center}

The first column contains the subgroups of $\Aut(A)$ grouped by conjugacy classes.
The next three columns contain the Galois group of the three extensions $K/k,K_2/k,K_1/k$.
The column labelled ``Type'' gives a name to each type of structure,
consisting of the Galois group of the field of definition of the Galois action
and an upper letter that distinguishes between the different types having the same group.
Finally, the last column gives some remarks on the fields
$K_1=k(\sqrt u)$ and $K_2=k(\sqrt v)$.

A $G_k$-group is completely determined, up to isomorphism,
by only giving its type and the field $K$ of definition of the action,
except for the two $V_4$-types,
where one has to specify which of the three quadratic subfields of $K$ plays the role of $K_2$
(three choices for every $K$),
and for type $D_8$,
where one must specify which of the two quadratic subfields different from $K_2$
plays the role of $K_1$
(two choices for every $K$).

For every row in the table, given extensions $K,K_1$ and $K_2$ of $k$
with the first containing the other two,
having the Galois groups that appear in that row,
and satisfying the conditions of the last column,
then there is a unique $G_k$-group structure on the group $D_8$
with the prescribed fieds as fields of definition.

\paragraph{Explicit construction of fields of definition.}
The following two lemmas introduce explicit constructions
of the fields of definition of $G_k$-modules isomorphic to $D_8$
in terms of radicals of elements of $k$:

\begin{lemma}\label{lemma_expression_field_D8}
Let $\car k\neq2$ and let $A$ be a $G_k$-group isomorphic to $D_8$.
Write $K_1=k(\sqrt u)$ and $K_2=k(\sqrt v)$ for elements $u,v\in k^*$.
Then, there exist elements $z\in k$ and $w\in k^*$ such that
$uv\equiv 1-z^2u\pmod{k^2}$ and with
\begin{equation}\label{expression_field_D8}
K=k(\sqrt v,\sqrt{\frac{w(1\pm z\sqrt u)}2}).
\end{equation}
Conversely, given elements $u,v,w\in k^*$ and $z\in k$ with $uv\equiv 1-z^2u\pmod{k^2}$
there exists a $G_k$-group isomorphic to $D_8$
with fields of definition $K_1=k(\sqrt u),K_2=k(\sqrt v)$ and
$K$ given by the expression (\ref{expression_field_D8}).
\end{lemma}

\begin{proof}
Assume a $G_k$-group $A$ isomorphic to $D_8$ is given.

Consider first the case $K_1=k$.
In this case the extension $K/k$ is trivial, quadratic or biquadratic,
and can be written as $K=k(\sqrt v,\sqrt m)$ for some $m\in k^*$.
Let $z\in k$ be an element such that $1-z^2u\equiv uv\pmod{k^2}$
(elements with this property always exist: every representation of a nonzero
square by the quadratic form $uvX^2+uZ^2$ produces such elements).
Choosing a square root $\sqrt u\in k$ with $1-z\sqrt u\neq0$
we may write the element $m$ as $m=w(1-z\sqrt u)/2$ for some $w\in k^*$.
From the condition that $1-z^2u\equiv uv\equiv v\pmod{k^2}$ it follows that
$w(1+z\sqrt u)/2\equiv mv\pmod{k^2}$ and since $K$ also contains $\sqrt{vm}$
we obtain that
$$K=k(\sqrt v,\sqrt m)=k(\sqrt v,\sqrt{\frac{w(1\pm z\sqrt u)}2})$$

Assume now that $K_1\neq k$.
Then if $[K:k]=4$ the field $K$ can be obtained by adjoining to $K_1$
the square root of a nonzero element and if $[K:k]=8$ then $\Gal(K/k)\simeq D_8$
and the field $K$ is the Galois closure of a quartic non-normal extension of $k$
obtained by adjoining to $K_1$ the square root of a nonzero element.
In either case let $\alpha=x+y\sqrt u$ be such a nonzero element.
After multiplying it by a square if necessary,
we may assume that $x$ is nonzero and then $\alpha$ may be written as $w(1-z\sqrt u)/2$
for some elements $w,z\in k$ with $w\neq0$.
Since the extension $K/k$ is normal,
the field $K$ also contains the square root of the conjugate element
$\overline\alpha=w(1+z\sqrt u)/2$ and then $K$ is obtained as (\ref{expression_field_D8}).
It only remains to be checked that $1-z^2u\equiv uv\pmod{k^{*2}}$ for the element $z$ we used,
but this congruence is a consequence of the fact that
the field $K$ contains the element
$\beta=\sqrt{\alpha}\sqrt{\overline\alpha}=\frac{w}{2}\sqrt{1-z^2 u}$
and an element of $G_k$ leaves $\beta$ fixed if and only if it fixes the element $\sqrt{uv}$.

For the converse, given elements $u,v,z,w$ with $1-z^2u\equiv uv\pmod{k^2}$ then
the field $K=k(\sqrt v,\sqrt{w(1\pm z\sqrt u)/2})$ is a normal extension of $k$
with Galois group isomorphic to a subgroup of $D_8$ and contains the fields
$K_1=k(\sqrt u)$ and $K_2=k(\sqrt v)$ as subfields.
Using Table 1, classifying the $G_k$-group structures on $D_8$,
it is easy to check that there is a structure with fields of definition the given fields.
\end{proof}

\begin{lemma}\label{lemma_different_w}
Let $\car k\neq2$ and let $A$ be a $G_k$-group isomorphic to $D_8$.
Let $u,v,w\in k^*$ and $z\in k$ elements determining the fields $K,K_1$ and $K_2$
as in the previous lemma.
Given an element $w'\in k^*$ there exists a $z'\in k$ such that $u,v,z',w'$ give
the same fields if, and only if,
$$(-v,w)=(-v,w')\qquad\text{in}\quad\Br_2(k)$$
\end{lemma}

\begin{proof}
Consider first the case $K_1=k$, that is equivalent to $u\in k^{*2}$.
For every element $z\in k$ with $1-z^2u\equiv uv\pmod{k^2}$ one has
the following identity
\begin{equation}\label{aqui}
(-v,\frac{1-z\sqrt u}2)=1
\end{equation}
in $\Br_2(k)$, provided that $1-z\sqrt u$ is nonzero.
This is a consequence of the following identity
$$1-\frac{1+z\sqrt u}2=\frac{1-z\sqrt u}2\qquad\Rightarrow\qquad
  \frac{1-z\sqrt u}2-\left(\frac{s}2\right)^2uv=\left(\frac{1-z\sqrt u}2\right)^2,$$
obtained by multiplying the self-evident identity on the left by $(1-z\sqrt u)/2$.
If the field $K$ admits the two expressions
$$K=k(\sqrt v,\sqrt{\frac{w(1-z\sqrt u)}2})=k(\sqrt v,\sqrt{\frac{w'(1-z'\sqrt u)}2})$$
then $w(1-z\sqrt u)$ and $w'(1-z'\sqrt u)$ differ in a square of $k$ or a square times $v$
and using twice the identity (\ref{aqui}) we have
$$(-v,w)=(-v,w\frac{1-z\sqrt u}2)=(-v,vw\frac{1-z\sqrt u}2)=
  (-v,w'\frac{1-z'\sqrt u}2)=(-v,w').$$
Conversely, if $(-v,w)=(-v,w')$ then using again (\ref{aqui}) we see that
$(-v,ww')=(-v,mw')=1$ with $m=w(1-z\sqrt u)/2$
and hence the quadratic form $mw'X^2+vmw'Y^2$ represents nonzero squares.
Let $x\in k$ be an element with $mw'x^2+vmw'y^2=1$.
If we write $mw'x^2=(1-z'\sqrt u)/2$ for some element $z'\in k$,
then $vmw'y^2=(1+z'\sqrt u)/2$ and it follows that $1-{z'}^2u\equiv v\equiv uv\pmod{k^2}$.
Moreover,
$$K=k(\sqrt v,\sqrt m)=k(\sqrt v,\sqrt{w'})=k(\sqrt v,\sqrt{\frac{w'(1-z'\sqrt u)}2}).$$
Let now $K_1\neq k$. Then
$$K=k(\sqrt v,\sqrt{\frac{w(1-z\sqrt u)}2}=k(\sqrt v,\sqrt{\frac{w'(1-z'\sqrt u)}2}$$
if, and only if, the two elements $w(1-z\sqrt u)/2$ and $w'(1-z'\sqrt u)/2$
differ in a square of $K_1$.
This condition is equivalent to the existence of elements $a,b\in k$ with
$$w(1-z\sqrt u)(a+b\sqrt u)^2=w'(1-z'\sqrt u)$$
for some $z'\in k$, and this is equivalent to the fact that the quadratic form
$$ww'(a^2+b^2u-2abzu)=ww'(a^2+u(b-az)^2-ua^2z^2)=ww'(1-z^2u)a^2+wu(b-az)^2$$
represents $1$ over $k$, which is equivalent to the identity
$$(ww'(1-z^2u),ww'u)=(ww'uv,ww'u)=(-v,ww'u)=(-v,ww')=1$$
that is equivalent to $(-v,w)=(-v,w')$ in $\Br_2(k)$.
\end{proof}

\paragraph{Realization of $G_k$-groups isomorphic to $D_8$ as groups of matrices.}
Since the groups of automorphisms of curves of genus 2 are $G_k$-groups that
can be represented as sub-$G_k$-groups of $\GL_2(\overline k)$,
we study these representations.

\begin{proposition}\label{proposition_matrices_D8}
Let $\car k\neq2$ and let $A$ be a sub-$G_k$-group of $\GL_2(\overline k)$ isomorphic to $D_8$.
For every $u,v\in k^*$ such that $K_1=k(\sqrt u)$ and $K_2=k(\sqrt v)$
there exists an element $z\in k$ such that
the group $A$ is $\GL_2(k)$-conjugated of the group generated by the two matrices
\begin{equation}\label{matrices_D8}
  U=\begin{pmatrix} \alpha & \beta \\ \beta/v & -\alpha \end{pmatrix},\qquad
  V=\begin{pmatrix} 0 & -\sqrt v \\ \frac{1}{\sqrt v} & 0 \end{pmatrix}
\end{equation}
with
\begin{equation}\label{matrices_D82}
\alpha=\sqrt{\frac{1-z\sqrt u}2},\qquad\beta=\sqrt{\frac{v(1+z\sqrt u)}2}.
\end{equation}
In that case the congruence $1-z^2u\equiv uv\pmod{k^2}$ is satisfied.

Conversely, given elements $u,v\in k^*$ and $z\in k$ with $1-z^2u\equiv uv\pmod{k^2}$
the matrices $U,V$ defined by (\ref{matrices_D8}) and (\ref{matrices_D82})
generate a group isomorphic to $D_8$ that is invariant by Galois action.
\end{proposition}

\begin{proof}
Let $A$ be a $G_k$-subgroup of $\GL_2(\overline k)$ isomorphic to $D_8$.
Since the central element $V^2$ of order 2 must be represented by the matrix $-1$,
the group is generated by matrices $U$ and $V$ with the relations
$$A=\langle\ U,V\ |\ U^2=1,V^2=-1,UV=-VU\ \rangle.$$
The matrix $V$ is defined over $k(\sqrt v)$ and
the matrix $V_0=\sqrt v V$ is fixed by every Galois automorphism
and has characteristic polynomial $X^2+v$.
After conjugation of the group $A$ by a matrix of $\GL_2(k)$ if necessary
we may assume that $V_0$ is the companion matrix of that polynomial
and then the matrix $V$ is
$$V=\begin{pmatrix} 0 & -\sqrt v \\ 1/\sqrt v & 0 \end{pmatrix}.$$

The matrices $U\in\GL_2(\overline k)$ of order two satisfying $UV=-VU$
are the matrices of the form
$$U=\begin{pmatrix} \alpha & \beta \\ \beta/v & -\alpha \end{pmatrix},\qquad
    \alpha,\beta\in\overline k,\qquad \alpha^2+\frac{\beta^2}{v}=1.$$
The matrix $U$ can have at most the four Galois conjugates $\pm U,\pm UV$.
Hence, it is defined over a field $K_U$ that contains field $K_1=k(\sqrt u)$
and such that the extension $K_U/K_1$ is either trivial or quadratic,
and the elements of $\Gal(K_U/K_1)$ send $U$ to $\pm U$.
It follows that $\alpha^2$ and $\beta^2$ belong to the field $K_1$.

If $K_1=k$ then we just define $z=(1-2\alpha^2)/\sqrt u$
and the elements $\alpha,\beta$ are expressed
in terms of this element by the formulas (\ref{matrices_D82}).

Let $K_1\neq k$.
Then every $\sigma\in G_k$ that is the nontrivial automorphism on $K_1$
acts on the matrix $U$ as ${}^\sigma U=\pm UV$
and since the matrix $UV$ is of the form
$$UV=\begin{pmatrix} \frac{\beta}{\sqrt v} & -\sqrt v\alpha \\
                     \frac{-\alpha}{\sqrt v} & -\frac{\beta}{\sqrt v} \end{pmatrix}$$
then ${}^\sigma\alpha^2=\beta^2/v$.
Now if we express $\alpha^2$ as an element of $k(\sqrt u)$,
say $\alpha^2=(x-z\sqrt u)/2$ with $x,z\in k$,
then $\beta^2$ is the element $v(x+z\sqrt u)/2$ of $k(\sqrt u)$.
From the identity $\alpha^2+\beta^2/v=1$ it follows that necessarily
$x$ must be $1$ and $\alpha$ and $\beta$ are given by the expressions
(\ref{matrices_D82}).

For every $\sigma\in G_k$ we have ${}^\sigma(\alpha\beta)=\pm\alpha\beta$
depending on the action of $\sigma$ on the field $K_1$.
Hence $\alpha\beta$ differs from $\sqrt u$ by an element of $k$ and
$$(\alpha\beta)^2=\frac{v(1-z^2u)}{4}\equiv u\pmod{k^{2}}\quad\Rightarrow\quad
  1-z^2u\equiv uv\pmod{k^{2}}.$$

As for the converse, let $u,v,z$ be elements of $k$ satisfying the requiered congruence.
Then it is immediate to check that the matrices $U$ and $V$
defined by (\ref{matrices_D8}) and (\ref{matrices_D82})
generate a group isomorphic to $D_8$.
If $\sigma\in G_k$ then ${}^\sigma V=\pm V$. Let
$$\overline\alpha=\sqrt{\frac{1+z\sqrt u}2}=\frac{\beta}{\sqrt v},\qquad
  \overline\beta=\sqrt{\frac{v(1-z\sqrt u)}2}=\alpha\sqrt v.$$
The Galois conjugates of $\alpha$ belong to the set $\{\pm\alpha,\pm\overline\alpha\}$.
If ${}^\sigma\alpha=\pm\alpha$ then ${}^\sigma\beta=\pm\beta$ and ${}^\sigma U=\pm U$;
if ${}^\sigma\alpha=\pm\overline\alpha$ then ${}^\sigma\beta=\mp\overline\beta$
and ${}^\sigma U=\pm UV$.
Hence the group $\langle U,V\rangle$ is closed by the Galois action.
\end{proof}

\begin{theorem}\label{theorem_obstruction_D8}
A $G_k$-group isomorphic to $D_8$ can be realized as a sub-$G_k$-group of $\GL_2(\overline k)$
if, and only if,
$$(-v,w)\qquad\text{is trivial in}\qquad\Br_2(k)$$
with $v,w\in k^*$ as in lemma \ref{lemma_expression_field_D8}.
\end{theorem}

\begin{proof}
From proposition \ref{proposition_matrices_D8} it follows that the $G_k$-group
is isomorphic to a $G_k$-group of matrices if, and only if, the field of definition
can be written as in lemma \ref{lemma_expression_field_D8} with the value $w=1$.
By lemma \ref{lemma_different_w}, given any elements $z,w$ determining the
field of definition of the corresponding action as in lemma \ref{lemma_expression_field_D8},
there exist elements $z',w'$ determining the same field with $w'=1$ if, and only if,
$(-v,w)=(-v,1)=1$.
\end{proof}

\paragraph{Galois structures on $D_{12}$.}
Now we consider the group $A\simeq D_{12}$, with presentation
$$A=\langle\ U,V\ |\ U^2=V^6=1,VU=UV^5\ \rangle.$$
The group $\Aut(A)$ is also isomorphic to $D_{12}$ and
it is generated by the two automorphisms $s$ and $t$ defined by
\begin{equation}\label{automorphisms_D12}
\begin{aligned} s(U) &= U  \\ s(V) &= V^5 \end{aligned}\qquad
\begin{aligned} t(U) &= UV \\ t(V) &= V \end{aligned}
\end{equation}
with the relations $s^2=t^6=1,ts=st^5$.
The inner automorphisms of $\Aut(A)$ are the subgroup of $\Aut(\Aut(A))$,
isomorphic to the dihedral group $D_6$,
generated by the conjugation by $s$ and the conjugation by $t$.

The characteristic subgroups of $A$ different form the trivial subgroup and from $A$ itself
are the subgroups $Z(A)=\langle V^3\rangle\simeq C_2$,
$\langle V^2\rangle\simeq C_3$ and $\langle V\rangle\simeq C_6$.
Every $G_k$-group structure on $A$ induces $G_k$-group structures on these groups
and also on the quotients $A'=A/Z(A)\simeq D_6$, $A/\langle V^2\rangle\simeq V_4$
and $A/\langle V\rangle\simeq C_2$. We remark that $A/\langle V^2\rangle=A^{\ab}$.

For every $G_k$-group structure on $A$ with field of definition $K$
we will denote by $K_1,K_2$ and $K_3$ the fields of definiton of the induced actions
on the groups $A/\langle V^2\rangle,\langle V\rangle$ and $A'$ respectively.
We remark that the field of definition of $\langle V\rangle$ is the same than
that of $\langle V^2\rangle$.
The Galois group of the extension $K/k$ is isomorphic to a subgroup of $D_{12}$,
the subfield $K_3/k$ has Galois group isomorphic to a subgroup of $D_6$
and the two extensions $K_i/k$ for $i=1,2$ can be either trivial or quadratic,
with $K_2\subseteq K_3$.

The following table classifies the $G_k$-group structures on $D_{12}$ by types,
in an analogous way than we did before for the group $D_8$

$$\begin{array}{|c|c|c|c|c|c|c|}
\hline
\im(\rho) & \Gal(K/k) & \Gal(K_3/k) & \Gal(K_2/k) & \Gal(K_1/k) & \text{Type} & \text{remarks} \\
\hline\hline
1                        &    I &   I &   I &   I & I & \\
\hline
\langle t^3 \rangle      &  C_2 &   I & I & C_2 & C_2^A & \\
\hline
\langle s\rangle,\langle st\rangle,\langle st^2\rangle &
   C_2 & C_2 & C_2 &  I & C_2^B & \\
\hline
\langle st^3 \rangle,\langle st^4\rangle,\langle st^5\rangle &
   C_2 & C_2 & C_2 & C_2 & C_2^C & K_1=K_2 \\
\hline
\langle t^2 \rangle      &  C_3 & C_3 &  I &   I & C_3 & \\
\hline
\langle t \rangle        &  C_6 &  C_3 & I & C_2 & C_6 & \\
\hline
\langle s,t^3 \rangle,\langle st,t^3\rangle,\langle st^2,t^3\rangle &
   V_4 & C_2 & C_2 & C_2 & V_4 & K_1\neq K_2 \\
\hline
\langle s,t^2 \rangle    &  D_6 & D_6 & C_2 &   I & D_6^A & \\
\hline
\langle st^3,t^2 \rangle &  D_6 & D_6 & C_2 & C_2 & D_6^B & K_1=K_2 \\
\hline
\langle s,t \rangle      & D_{12} & D_6 & C_2 & C_2 & D_{12} & K_1\neq K_2\\
\hline
\end{array}$$
\begin{center}{\bf Table 2.} $G_k$-group structures on $D_{12}$\end{center}
\medskip

A $G_k$-group isomorphic to $D_{12}$ is completely determined,
up to isomorphism, by giving its type
and the field of definition of the action $K/k$,
except for type $V_4$, where one needs to specify
which of the three quadratic subfields of $K$ are the fields $K_1$ and $K_2$
(six choices for a given $K$),
and for type $D_{12}$, where one must specify which of the two quadratic subfields
different from $K_2$ plays the role of $K_1$ (two choices for a given $K$).

Conversely for every row in the table,
given extensions $K,K_3,K_2,K_1$ of $K$ with the first containing the other three,
$K_3$ containing $K_2$,
with Galois groups as a in that row and satisfying the restrictions of the last column,
there is a unique $G_k$-group structure on the group $D_{12}$
with the given fields as fields of definition of the corresponding induced actions.

\paragraph{Realization of $G_k$-groups isomorphic to $D_{12}$ as groups of matrices.}
Now we study the representations of these $G_k$-groups as groups of matrices
in an analogous manner that we did before for the case $D_8$:

\begin{proposition}\label{proposition_matrices_D12}
Let $\car k\neq2,3$ and let $A$ be a sub-$G_k$-group of $\GL_2(\overline k)$
isomorphic to $D_{12}$.
For every element $u\in k^*$ such that $K_1=k(\sqrt u)$ and $K_2=k(\sqrt v)$
there exists an element $z\in k$ such that
the group $A$ is $\GL_2(k)$-conjugated of the group generated by the two matrices
\begin{equation}\label{matrices_D12}
U=\frac{1}{\sqrt u}\begin{pmatrix} \alpha & \beta \\ \beta/v & -\alpha \end{pmatrix},
\qquad
V=\frac12\begin{pmatrix} 1 & \sqrt v \\ -\frac3{\sqrt v} & 1 \end{pmatrix},
\end{equation}
with
\begin{equation}\label{matrices_D122}
\alpha\text{\ \ a root of}\quad X^3-\frac{3u}4X-\frac{z}4,
\qquad\text{and}\qquad
\beta\text{\ \ a root of}\quad\alpha^2+\frac{\beta^2}v=u.
\end{equation}
In that case the congruence $u^3-z^2\equiv3v\pmod{k^2}$ is satisfied.

Conversely, given elements $u,v\in k^*$ and $z\in k$ with $u^3-z^2\equiv3v\pmod{k^2}$
the matrices $U,V$ defined by (\ref{matrices_D12}) and (\ref{matrices_D122})
generate a group isomorphic to $D_{12}$ invariant by Galois action.
\end{proposition}

\begin{proof}
Let $A$ be a $G_k$-subgroup of $\GL_2(\overline k)$ isomorphic to $D_{12}$.
Since the central element $V^3$ of order $2$ must be represented by the matrix $-1$
the group can be written as
$$A=\langle\ U,V\ |\ U^2=1,V^3=-1,VU=-UV^2\ \rangle$$
The matrix $V$ is defined over the field $K_2$ and the Galois action either
leaves it unchanged, if $K_2=k$, or interchanges $V$ with $V^5$, if $K_2\neq k$.
Since the characteristic polynomial of $V$ is $X^2-X+1$
it follows that $V$ must be $\GL_2(k)$-conjugated of the matrix
$$\frac12\begin{pmatrix} 1 & \sqrt v \\ -\frac3{\sqrt v} & 1 \end{pmatrix}.$$

The matrices $U\in\GL_2(\overline k)$ of order two satisfying $VU=-UV^2$
are the matrices of the form
$$U=\frac{1}{\sqrt u}\begin{pmatrix} \alpha & \beta \\ 3\beta/v & -\alpha \end{pmatrix},
\qquad\alpha,\beta\in\overline k,\qquad\alpha^2+\frac{3\beta^2}{v}=u.$$
Since the set of matrices $\{\sqrt u\ U,\sqrt u\ UV^2,\sqrt u\ UV^4\}$
is closed by Galois action,
it follows that the upper left entries of these three matrices are the roots of a polynomial
with coefficients in $k$.
The identity $\alpha^2+3\beta^2/v=u$ implies that this polynomial must be of the form
$$X^3-\frac{3u}4X-\frac{z}4$$
for some element $z\in k$, and $\alpha$ and $\beta$
satisfy the conditions (\ref{matrices_D122}).

The congruence $u^3-z^2\equiv3v\pmod{k^2}$ is a consequence of the fact that the roots of
the polynomial $X^3-\frac{3u}4X-\frac{z}4$ generate the field $K_3$
whose unique quadratic subfield is $K_2$ and the discriminant of the polynomial
is $3^32^{-4}(u^3-z^2)$.

Conversely, let $u,v,z$ be elements satisfying the congruence,
and $\alpha,\beta\in\overline k$ be elements satisfying the conditions (\ref{matrices_D122}).
Then it is immediate to check that the matrices $U$ and $V$ defined by (\ref{matrices_D12})
generate a group isomorphic to $D_{12}$.
If $\sigma\in G_k$ then ${}^\sigma V\in\{V,V^2\}$.
If the polynomial $X^3-\frac{3u}4X-\frac{z}4$ has multiple roots then
$\alpha$ and $\beta$ are in an at most quadratic extension of $k$ and it is clear that
$\langle U,V\rangle$ is closed by Galois action.
If that polynomial is separable then $\beta$ can be written in terms of $\alpha$ as
$$\beta=\frac{u^2+z\alpha-2u\alpha^2}{3s}$$
for an element $s\in k^*$ such that $u^3-z^2=3vs^2$.
In terms of this element, the roots of the polynomial different from $\alpha$ are
$$\alpha',\alpha''=\frac{u^2+z\alpha-2u\alpha^2\pm s\alpha\sqrt v}{\pm\sqrt v}$$
and one has ${}^\sigma U\in\{U,UV^2,UV^4\}$
depending on ${}^\sigma\alpha\in\{\alpha,\alpha',\alpha''\}$.
Hence the group $\langle U,V\rangle$ is closed by the Galois action.
\end{proof}

\begin{theorem}\label{theorem_obstruction_D12}
A $G_k$-group isomorphic to $D_{12}$ can be realized as a sub-$G_k$-group of $\GL_2(\overline k)$
if, and only if,
$$(u,-3v)\qquad\text{is trivial in}\qquad\Br_2(k)$$
with $u$ and $v$ elements such that $K_1=k(\sqrt u)$ and $K_2=k(\sqrt v)$.
\end{theorem}

\begin{proof}
If the group can be realized then by proposition (\ref{proposition_matrices_D12})
there exists an element $z\in k$ satisfying the congruence $u^3-z^2\equiv3v\pmod{k^2}$.
Then, the quadratic form $uX^2-3vY^2$ represents a square over $k$
and this implies the identity $(u,-3v)=1$.

Let a $G_k$-group isomorphic to $D_{12}$ be given
and assume that the condition $(u,-3v)=1$ is satisfied.

If $3\nmid[K:k]$ then $K=k(\sqrt u,\sqrt v)$.
The identity on the Brauer group
implies that the quadratic form $uX^2-3vY^2$ represents squares over $k$,
and from a representation with nonzero $X$ we obtain elements $\alpha,\beta\in k$
satisfying $\alpha^2+3\beta^2/v=u$.
Then the matrices defined by formulas (\ref{matrices_D12}) using these elements
$\alpha,\beta$ generate a group of matrices isomorphic to the given $G_k$-group.

Now consider the case $3\mid[K:k]$. Then $K_3/k$ is cyclic of degree 3 or dihedral of degree 6.
We want to see that this extension is the field of decomposition of some polynomial
of the form $X^3-\frac{3u}4X-\frac{z}4$.
Let $f(X)=X^3-aX-b\in k[X]$ be any irreducible polynomial with $K_3$ as field of decomposition,
and let $\alpha\in\overline k$ be a root of this polynomial.
Then $4a^3-27b^2\equiv v\pmod{k^{*2}}$ since the field $K_2=k(\sqrt v)$ is contained in $K_3$.
If $a=0$ then the extension $k(\alpha)$ can always be written as the field of decomposition
of a polynomial as claimed.
Assume $a\neq0$.
Computing the irreducible polynomial of a generic element of $k(\alpha)$ of zero trace
one obtains the polynomials $f_1(X)=X^3-a_1X-b_1$ with $a_1=ax^2+3avy^2$ and $x,y\in k$.
The formula for the discriminant of $f(X)$ in terms of $a$ and $b$
shows that the quadratic form $avX^2-3vY^2$ represents squares in $k$ and hence $(av,-3v)=1$.
Since by hypothesis also $(u,-3v)=1$ it follows that $(avu,-3v)=(3au,-3v)=1$
and from this identity one deduces the existence of elements $x,y\in k$
with $a_1=ax^2+3avy^2=3u/4$. Taking $z=4b_1$ for the corresponding coefficient $b_1$
the extension $K_3/k$ is the field of decomposition of a polynomial
of the type $X^3-\frac{3u}4X-\frac{z}4$ as claimed.
Now, the matrices defined by formulas (\ref{matrices_D12}) using a root $\alpha$
of this polynomial and a root $\beta$ of $\alpha^2+\beta^2/v=u$ generate a
sub-$G_k$-module of $\GL_2(\overline k)$ isomorphic to the module we started with.
\end{proof}

\paragraph{Two technical lemmas.}
We end this section with two technical lemmas.
The first one studies a particular cohomology class with values in a $G_k$-group
isomorphic to $D_8$ that will be used later because it produces the
unique non-hyperelliptic twists that maintain the Galois action on the group of automorphisms.
The second lemma shows that two isomorphic $G_k$-groups realized as groups of matrices
that are isomorphic to $D_8$ or $D_{12}$
are conjugated by some matrix with coefficients in the base field $k$.

\begin{lemma}\label{lemma_special_cocycle}
Let $A$ be a $G_k$-module isomorphic to $D_8$ as a group
and let $K_2=k(\sqrt v)$ for some $v\in k^*$.
The 1-cocycle $\Xi$ defined by
$$\Xi_\sigma=\begin{cases}
1,\quad {}^\sigma\sqrt v=\sqrt v,\\
V,\quad {}^\sigma\sqrt v=-\sqrt v.\end{cases}$$
is cohomologous to a one-cocycle with values in $\{\pm1\}$
if, and only if, $K_2=k$ or $K_2=K_1$.
\end{lemma}

\begin{proof}
We just multiply the cocycle $\Xi$ by all the coboundaries,
which are the maps of the form $\sigma\mapsto {}^\sigma W\cdot W^{-1}$ for elements $W\in A$,
and check that one gets a cocycle with values in $\{\pm1\}$ for some $W$
exactly when one of the conditions $K_2=k$ or $K_2=K_1$ hold.
\end{proof}

\begin{lemma}\label{lemma_isomorphic_Gkgroups}
Two $G_k$-subgroups of $\GL_2(\overline k)$ isomorphic to $D_8$
that are isomorphic as $G_k$-groups are $\GL_2(k)$-conjugate.
The same holds if they are both isomorphic to $D_{12}$.
\end{lemma}

\begin{proof}
Let $A_1$ and $A_2$ be subgroups of $\GL_2(\overline k)$
that are isomorphic as $G_k$-groups, and are isomorphic to $D_8$ as groups.

Since $D_8$ has a unique faithful two-dimensional representation,
every group isomorphism between $A_1$ and $A_2$
is obtained by conjugation by some matrix of $\GL_2(\overline k)$.
Let $f:A_1\to A_2$ be an isomorphism as $G_k$-groups,
and let $M\in\GL_2(\overline k)$ be a matrix such that $f(W)=MWM^{-1}$ for all $W\in A_1$.
Then, and since $f$ respects the Galois action, we have
$$M{}^\sigma WM^{-1}=f({}^\sigma W)={}^\sigma f(W)={}^\sigma M{}^\sigma W{}^\sigma M^{-1}$$
for every $W\in A_1$ and $\sigma\in G_k$.

For every $\sigma\in G_k$,
the conjugation by the matrix $M^{-1}\cdot{}^\sigma M$
acts as the identity on the group $A_1$ and hence is an homotety,
which we identify with an element of $\overline k^*$.
The map $\sigma\mapsto M^{-1}\cdot{}^\sigma M\in k^*$ is a 1-cocycle of $G_k$
with values in the multiplicative group $\overline k^*$ and by the theorem 90 of Hilbert
its cohomology class is trivial.
Let $a\in\overline k^*$ be an element such that $M^{-1}\cdot{}^\sigma M=a\cdot{}^\sigma a^{-1}$.
Then the matrix $Ma$ is Galois invariant, hence it belongs to $\GL_2(k)$,
and the isomorphism $f$ is also obtained by conjugation by this matrix.

The same proof works for the case of groups isomorphic to $D_{12}$ and, in fact,
for every group that has a unique faithful 2-dimensional representation
with centralizer in $\GL_2(\overline k)$ consisting only of homoteties.
\end{proof}

\section{Classification of twists}

The main object of this section is the classification of the $k$-twists
of a curve of genus 2 defined over a field $k$ with group of automorphisms isomorphic to
$D_8$ or $D_{12}$.
Our strategy will be to start with the classification of hyperelliptic twists,
corresponding to the group $H^1(G_k,\{\pm1\})/\delta({A'}^{G_k})$
on the left of the exact sequence (\ref{cohomology_exact_sequence}), which is the easy part,
and then study general twists modulo hyperelliptic twists,
that in some sense can be seen as the classification of the elements of the
set $H^1(G_k,A')[\Delta]$ on the right of the same exact sequence.

The computation of the group $H^1(G_k,\{\pm1\})/\delta({A'}^{G_k})$
is reduced to that of the finite image $\delta({A'}^{G_k})$
for every possible $G_k$-structure on $A$, and the result is given in the following

\begin{proposition}\label{proposition_hyperelliptic_twists}
Let $\car k\neq2$ and
let $A$ be a sub-$G_k$-group of $\GL_2(\overline k)$ isomorphic to $D_8$. Then,
$$H^1(G_k,\{\pm1\})/\delta({A'}^{G_k})\simeq\begin{cases}
  k^*/(k^{*2}\cdot vk^{*2}\cdot mk^{*2}\cdot vmk^{*2}),&
  \quad\text{if}\quad K_1=k,\\
  k^*/(k^{*2}\cdot vk^{*2}),&
  \quad\text{if}\quad K_1\neq k,\end{cases}$$
with $K_2=k(\sqrt v)$ and, if $K_1=k$, then $K=k(\sqrt v,\sqrt{m})$, for elements $v,m\in k^*$.

Let $\car k\neq2,3$ and let $A$ be a sub-$G_k$-group of $\GL_2(\overline k)$
isomorphic to $D_{12}$. Then,
$$H^1(G_k,\{\pm1\})/\delta({A'}^{G_k})\simeq\begin{cases}
  k^*/(k^{*2}\cdot uk^{*2}),&
  \quad\text{if}\quad 3\nmid[K:k],\\
  k^*/k^{*2},&
  \quad\text{if}\quad 3\mid[K:k],\end{cases}$$
with $K_1=k(\sqrt u)$ for an element $u\in k^*$.
\end{proposition}

Now we define two familes of curves by giving explicit hyperelliptic equations for them,
depending on parameters of $k$ that satisfy one relation:

\begin{definition}\label{definition_rational_models}
Let $\car k\neq 2$ and also $\neq3$ in the $D_{12}$ case.
Let $u,v\in k^*$ and $z,s\in k$ be elements satisfying the identity
$$1-z^2 u=s^2uv\qquad\quad{or}\qquad u^3-z^2=3s^2v$$
in the $D_8$ case and the $D_{12}$ case respectively. We define the curves:
$$\begin{aligned}
C^{(8)}_{u,v,z,\pm}:\quad Y^2 =
&\ (1\pm 2uz) X^6 \mp 8suv X^5 + v(3\mp 10uz) X^4 \\
&\ + v^2(3\mp 10uz) X^2 \pm 8suv^3 X + v^3(1\pm 2uz)\\
\phantom{a} & \\
C^{(12)}_{u,v,z}:\quad Y^2 =
&\ 27(u+2z) X^6 - 324sv X^5 + 27v(u-10z) X^4 + 360sv^2 X^3 \\
&\ + 9v^2(u+10z) X^2 - 36sv^3 X + v^3(u-2z)\end{aligned}$$
\end{definition}

We remark that $s$ is determined by $u,v,z$ up to a sign,
and that a change of sign produces $k$-isomorphic curves.
In fact, when $s\neq0$ the curves obtained from different nonzero values of $s$
are $k$-isomorphic to each other,
hence we can safely replace $s$ by $1$ in the equation
without changing the isomorphism class.

If $u\neq0,-1/4$ then these equations correspond to curves of genus 2
of the type we are interested in, since

\begin{proposition}
The curves $C^{(8)}_{u,v,z,\pm}$ have group of automorphisms isomorphic to $D_8$
generated by the matrices $U,V$ of proposition \ref{proposition_matrices_D8},
and their absolute invariant is $t(C^{(8)}_{u,v,z,\pm})=u$.

The curves $C^{(12)}_{u,v,z}$ have group of automorphisms isomorphic to $D_{12}$
generated by the matrices $U,V$ of proposition \ref{proposition_matrices_D12},
and their absolute invariant is $t(C^{(12)}_{u,v,z})=u$.
\end{proposition}

\begin{proof}
Direct computation:
one checks that the substitutions corresponding to the matrices $U$ and $V$
are automorphisms of the given equations,
and computes the absolute invariant $t$ in terms of invariants of sectic forms
by the formulas given in propositions \ref{proposition_absolute_invariant_D8}
and \ref{proposition_absolute_invariant_D12}.
\end{proof}

\paragraph{Classification up to hyperelliptic twist in the $D_8$ case.}
In the $D_8$ case we defined two curves for every set of parameters $u,v,z$ of $k$,
depending on a coherent choice of signs in the equation defining $C^{(8)}_{u,v,z,\pm}$.
The effect of the change of sign is:

\begin{proposition}\label{proposition_change_sign_equation}
The two curves with equations $C^{(8)}_{u,v,z,+}$
and $C^{(8)}_{u,v,z,-}$ are
twists of each other by the element $\Xi\in H^1(G_k,A)$ of lemma \ref{lemma_special_cocycle},
up to hyperelliptic twist.
\end{proposition}

\begin{proof}
One checks that the matrix
$$\Phi=\frac{1}{\sqrt 2}\begin{pmatrix} 1 & \sqrt v \\ \frac{-1}{\sqrt v} & 1 \end{pmatrix}
  \in\GL_2(\overline k)$$
gives an isomorphism between these two curves.
The map $\sigma\mapsto{}^\sigma\Phi\cdot\Phi^{-1}$ is the product of the one-cocycle $\Xi$
of lemma \ref{lemma_special_cocycle}
by the one-cocycle with values in $\{\pm1\}$ corresponding to the field $k(\sqrt2)$.
\end{proof}

Hence, the two curves are an hyperelliptic twist of each other if $K_2=k$ or $K_2=K_1k$,
and a non-hyperelliptic twist of each other if $k\neq K_2\neq K_1$.
The last possibility corresponds to Galois actions on $\Aut(A)$
of types $C_2^B, V_4^A$ and $D_8$.

\begin{theorem}\label{main_theorem_D8}
Let $C_i/k, i=1,2$ be curves of genus 2 with $A_i=\Aut(C_i)\simeq D_8$.
If the two curves $C_i$ are $\overline k$-isomorphic to each other
and the two groups $A_i$ are isomorphic to each other as $G_k$-groups
then the curves $C_i$ differ by a twist of the type $\xi$ or $\xi\Xi$,
with $\xi$ an hyperelliptic twist
and $\Xi$ the twist of lemma (\ref{lemma_special_cocycle}).
\end{theorem}

\begin{proof}
By lemma \ref{lemma_isomorphic_Gkgroups} we may assume,
changing one of the curves by a $k$-isomorphic one if necessary,
that the groups of automorphisms of the two curves,
viewed as subgroups of $\GL_2(\overline k)$, are the same group.
By proposition \ref{proposition_matrices_D8}, and again up to a $k$-isomorphism,
we may assume that this common group $A=\Aut(C_i)$
is a group generated by two matrices $U$ and $V$
as in proposition \ref{proposition_matrices_D8}.

Let $M\in\GL_2(\overline k)$ be
the matrix associated to an isomorphism between these two curves.
Conjugation by $M$ gives an automorphism of the group $A$.
We recall that the identity and the automorphism $t$ defined in (\ref{automorphisms_D8})
represent the two cosets of $\Aut(A)$ modulo inner automorphisms.
Up to multiply $M$ by an element of $A$, which gives another isomorphism between the curves,
we may assume that the conjugation by $M$ is either the identity or the automorphism $t$.

If conjugation by $M$ is the identity on $A$, then $M$ must be a scalar matrix.
For every $\sigma\in G_k$ the matrix $\xi_\sigma={}^\sigma M\cdot M^{-1}\in A$
is also scalar and hence it must be $\pm1$.
In this case, the two curves are hyperelliptic twists of each other by $\xi$.

Assume now that conjugation by $M$ gives the automorphism $t$.
Let $C_2'$ be the curve obtained applying to $C_2$ the transformation given
by the matrix $\Phi$ in the proof of proposition \ref{proposition_change_sign_equation};
the curves $C_2'$ and $C_2$ are twists corresponding to the element $\xi_2\Xi$,
with $\xi_2:G_k\to\{\pm1\}$ the one-cocycle corresponding to the extension $k(\sqrt2)$.
One checks that conjugation by $\Phi$ induces the automorphism $t$
of the group of matrices $\Aut(C_i)=\langle U,V\rangle$.
Then the curves $C_1$ and $C_2'$ also have the same group of automorphisms
and there is an isomorphism between them that induces by conjugation
the identity on $\langle U,V\rangle$.
By the previous argument the curves $C_2$ and $C_2'$ differ by an hyperelliptic twist $\xi$
and hence $C_1$ and $C_2$ are related by the twist $\xi\xi_2\Xi$.
\end{proof}

\begin{corollary}
Let $\car k\neq2$.
\begin{enumerate}
\item Every curve $C/k$ of genus 2 with $\Aut(C)\simeq D_8$ is isomorphic to
   an hyperelliptic twist of a curve $C^{(8)}_{u,v,z,\pm}$
   for some paremeters $u,v,z$ and sign $\pm$.
\item The pair $C^{(8)}_{u,v,z,\pm}$ is an hyperelliptic twist
   of the pair $C^{(8)}_{u',v',z',\pm}$ if, and only if,
   $$u=u',\quad k(\sqrt v)=k(\sqrt{v'})\quad\text{and}\quad
  k(\sqrt v,\sqrt{\frac{1\pm z\sqrt u}2})=k(\sqrt{v'},\sqrt{\frac{1\pm z'\sqrt{u'}}2})$$
\end{enumerate}
\end{corollary}

\begin{proof}
1. Given such a curve $C/k$ let $t$ be its absolute invariant, and let $A=\Aut(C)$.
Then, the field $K_1$ corresponding to the Galois action on $A$ is $k(\sqrt t)$.
Indeed, the curve of genus 2 with equation $C_t:Y^2=X^5+X^3+t\,X$
is isomorphic to the curve $C$,
and hence the $G_k$-group $A_t=\Aut(C_t)$
is twisted of $A$ by some one-cocycle $\xi:G_k\to A$.
Hence the Galois action on the quotients $A^{\ab}$ and $A_t^{\ab}$ are the same.
If $K_1$ denotes the field of definition of $A^{\ab}$ then $K_1=k(\sqrt t)$,
since $k(\sqrt t)$ is the field of definition of $A_t^{\ab}$.

Applying proposition \ref{proposition_matrices_D8} to the group $A$,
viewed as a group of matrices,
with the choice of parameters $u=t$ and $v\in k^*$ any element with $K_2=k(\sqrt v)$
we obtain an element $z\in k$ and an $s\in k$ such that $1-z^2t=s^2tv$.
The curves $C^{(8)}_{t,v,z,\pm}$ defined from these parameters
are isomorphic to the curve $C$, since they have the same absolute invariant,
and have group of automorphisms isomorphic to $A$ as $G_k$-groups,
since both groups are conjugated by means of a matrix of $\GL_2(k)$.
By theorem \ref{main_theorem_D8} and proposition \ref{proposition_change_sign_equation}
at least one of the two curves $C^{(8)}_{t,v,z,+}$ and $C^{(8)}_{t,v,z,-}$
differ from $C$ by an hyperelliptic twist.

2. The second part of the statement
is also an immediate consequence of theorem \ref{main_theorem_D8}
since the curves have absolute invariants $u$ and $u'$ respectively,
and the Galois action on their groups of automorphisms
is determined by the fields $K_1,K_2$ and $K$ corresponding to these actions.
\end{proof}

\paragraph{Classification up to hyperelliptic twist in the $D_{12}$ case.}

\begin{theorem}
Let $C_i/k, i=1,2$ be curves of genus 2 with $A_i=\Aut(C_i)\simeq D_8$.
If the two curves $C_i$ are $\overline k$-isomorphic to each other
and the two groups $A_i$ are isomorphic to each other as $G_k$-groups
then the curves $C_i$ differ by a hyperelliptic twist.
\end{theorem}

\begin{proof}
By lemma \ref{lemma_isomorphic_Gkgroups} we may assume,
changing one of the curves by a $k$-isomorphic one if necessary,
that the groups of automorphisms of the two curves,
viewed as subgroups of $\GL_2(\overline k)$, are the same group.
Again up to $k$-isomorphism we may assume that this common group $A=\Aut(C_i)$
is a group generated by two matrices $U$ and $V$
as in proposition \ref{proposition_matrices_D12}.

Let $M\in\GL_2(\overline k)$ be
the matrix associated to an isomorphism between these two curves.
Conjugation by $M$ gives an automorphism of the group $A$.
We recall that the identity and the automorphism $t^3$ defined in (\ref{automorphisms_D12})
represent the two cosets of $\Aut(A)$ modulo inner automorphisms.
Up to multiply $M$ by an element of $A$, which gives another isomorphism between the curves,
we may assume that the conjugation by $M$ is the identity or the automorphism $t^3$.

Let $\Phi$ be the matrix
$$\Phi=\begin{pmatrix} 0 & -v/3 \\ 1 & 0 \end{pmatrix}\in\GL_2(k).$$
Conjugation by the matrix $A$ produces also the automorphism $t^3$ of the group $A$.
Changing one of the curves by a $k$-isomorphic curve with isomorphism given by $\Phi$,
if conjugation by $M$ is the automorphism $t^3$ of $A$,
we may assume that conjugation by $M$ is the identity on $A$.
Then $M$ must be a scalar matrix.
For every $\sigma\in G_k$ the element $\xi_\sigma={}^\sigma M\cdot M^{-1}\in A$
is also scalar and hence it must be $\pm1$.
Then he two curves are hyperelliptic twists of each other by $\xi$.
\end{proof}

\begin{corollary}
Let $\car k\neq2,3$.
\begin{enumerate}
\item Every curve $C/k$ of genus 2 with $\Aut(C)\simeq D_{12}$ is isomorphic to
   an hyperelliptic twist of a curve $C^{(12)}_{u,v,z}$
   for some paremeters $u,v,z$,
\item The curve $C^{(12)}_{u,v,z}$ is an hyperelliptic twist
   of the curve $C^{(12)}_{u',v',z'}$ if, and only if,
$$u=u',\quad k(\sqrt v)=k(\sqrt{v'})\quad\text{and}\quad
  k(\sqrt v,\alpha)=k(\sqrt{v'},\alpha')$$
with $\alpha$ and $\alpha'$ roots of the polynomials
$$X^3-\frac{3u}4 X-\frac{z}4,\qquad X^3-\frac{3u'}4 X-\frac{z'}4.$$
\end{enumerate}
\end{corollary}

\begin{proof}
The proof is analogous to that of the corresponding corollary for the $D_8$ case.
Perhaps the only thing that should be remarked is that, also in this case,
the field $K_1$ corresponding to a $G_k$-group $A$ isomorphic to $D_{12}$
is the field of definition of the action on the quotient group $A^{\ab}$,
which is invariant by twisting the action,
and hence is an invariant of the $\overline k$-isomorphism class of the curve.
\end{proof}

\section{Elliptic quotients and abelian varieties of $\GL_2$-type}

We begin by recalling some definitions and facts
about $\Q$-curves and abelian varieties of $\GL_2$-type;
for details the reader may consult \cite{ribet} or \cite{quer}.
A $\Q$-curve is an elliptic curve defined over a number field
that is isogenous to all its Galois conjugates;
it is said to be completely defined over a number field $k$
if all the conjugates of the curve and the isogenies between them are defined over $k$.
An abelian variety of $\GL_2$-type is an abelian variety defined over $\Q$
such that the $\Q$-algebra $\Q\otimes\End_\Q(A)$ of endomorphisms up to isogeny
defined over $\Q$ is a number field of degree equal to $\dim A$.
The main reason for the interest in these varieties are the modularity conjectures:
\begin{itemize}
\item The abelian varieties of $\GL_2$-type are, up to $\Q$-isogeny,
   the $\Q$-simple factors of modular jacobians $J_1(N)$.
\item The $\Q$-curves are, up to isogeny,
   the one-dimensional factors of modular jacobians $J_1(N)$.
\end{itemize}
Ribet proved that the first implies the second
and that both conjectures would follow from Serre's conjecture
on 2-dimensional mod $p$ Galois representations.
The modularity of all the $\Q$-curves (isogenous to curves) defined over $\Q$
is now known from results by Wiles, Taylor-Wiles, and Breuil-Conrad-Diamond-Taylor.
Moreover, using Wiles' techniques, the modularity of some other families of $\Q$-curves
has been proved by Ellenberg-Skinner, Hasegawa-Hashimoto-Momose, and Hida.

In this section we show that the elliptic quotients
of the curves $C/\Q$ of genus 2 with $\End(C)\simeq D_8$ and $D_{12}$
are the $\Q$-curves of degrees 2 and 3, respectively.
Then, we characterize the curves $C/\Q$ of genus 2 in the families studied in this paper
with jacobian $J_C$ of $\GL_2$-type.
We remark that, assuming the modularity conjecture for $\Q$-curves of degrees 2 and 3,
the fact that the jacobian is of $\GL_2$-type is equivalent to have an $L$-series
that is a product of the $L$-series of two modular forms
for congruence subgroups $\Gamma_1(N)$.
Finally we show that several families of $\Q$-curves completely defined over a field $K$
can be covered by curves of genus 2 with covering also defined over $K$,
and see the relation between this fact and the property of the jacobian being
of $\GL_2$-type.

Although most facts on elliptic quotients of curves of genus 2 and on
endomorphisms of their jacobians given in this section could be stated over
any field $k$ with the usual restriction on the characteristic of previous sections,
we will work only for $k=\Q$ since the applications to modularity are for this field.

\begin{proposition}
Let $C/\Q$ be a curve of genus 2 with $A=\Aut(C)\simeq D_8$ or $D_{12}$
and let $K$ be the field of definition of $A$. Then
\begin{enumerate}
\item The curve $C$ has, up to isogeny, a unique elliptic quotient,
   which is a $\Q$-curve of degree 2 or 3, respectively.
\item For every non-hyperelliptic involution $W\in\Aut(C)$
   the quotient $E_W=C/\langle W\rangle$ is an elliptic curve
   and the morphism $C\to E_W$ is defined over the field of definition of $W$.
\item Two quotients $E_{W_1}$ and $E_{W_2}$ are isomorphic if, and only if,
   the involutions $W_1$ and $W_2$ are conjugated in $A$;
   in that case the isomorphism is defined over $K$.
\item The elliptic quotients $E_W$ with no complex multiplication are $\Q$-curves
   completely defined over the field $K$.
\end{enumerate}
\end{proposition}

\begin{proof}
See \cite[Section 2]{esquela} for a proof of these facts,
and also for an explicit construction of Weierstrass equations for the elliptic quotients
in terms of an hyperelliptic equation for the curve $C$ and the matrix corresponding
to the non-hyperelliptic involution $W$.
\end{proof}

If $\Aut(C)\simeq D_8$ there are four elliptic quotients
corresponding to the four hyperelliptic involutions $U,-U,UV,-UV$,
that belong to two classes of isomorphism.
When $\Aut(C)\simeq D_{12}$ the situation is analogous:
the six quotients corresponding to the involutions $U,V,UV^2,-U,-V,-UV^2$,
are grouped into two classes of isomorphism,
corresponding to the two classes of conjugation.
Then, in either case $\Aut(C)\simeq D_8$ or $\Aut(C)\simeq D_{12}$
there are only two elliptic quotients by non-hyperelliptic involutions up to isomorphism.
The following diagram shows the graph of the elliptic quotients,
whose edges represent isomorphisms (degree 1) and isogenies of degree 2 or 3
linking the curves, for the cases $D_8$ and $D_{12}$ respectively
$$
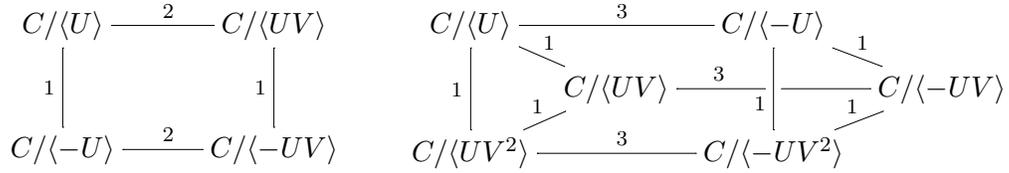
\xymatrix@+6pt{
  C/\langle U\rangle  \ar@{-}[r]^2 \ar@{-}[d]_1  & C/\langle UV\rangle  \ar@{-}[d]_1     \\
  C/\langle -U\rangle  \ar@{-}[r]^2            & C/\langle -UV\rangle                }
\qquad
\xymatrix@-17pt{
  C/\langle U\rangle  \ar@{-}[rr]^3\ar@{-}[dd]_1\ar@{-}[dr]^1 & & C/\langle -U\rangle
  \ar@{-}[dd]_{\sr{}{\sr{}{\sr{}{\sr{}1}}}}\ar@{-}[dr]^{1}              \\
  &      C/\langle UV\rangle  \ar@{-}'[r]^{\ \ \ \ \ 3}[rr]\ar@{-}[dl]_1 & & C/\langle -UV\rangle  \ar@{
-}[dl]_1\\
  C/\langle UV^2\rangle  \ar@{-}[rr]^3              & & C/\langle
  -UV^2\rangle
}$$
Let $t$ be the absolute invariant of the curve $C$,
and let $E_1$ and $E_2$ be two elliptic quotients representing the two classes of isomorphism.
One computes their $j$-invariants with the formula in \cite[Lemma 2.2]{esquela} and obtains
the following values
\begin{equation}\label{j_invariants_quotients}
(j_{E_1},j_{E_2})=\begin{cases}
\displaystyle\frac{2^6(3\mp10\sqrt t)^3}{(1\mp2\sqrt t)(1\pm2\sqrt t)^2},&\qquad \Aut(C)=D_8,\\
\phantom{\ }\\
\displaystyle\frac{2^8 3^3(2\mp5\sqrt t)^3(\pm\sqrt t)}{(1\mp2\sqrt t)(1\pm2\sqrt t)^3},&
\qquad \Aut(C)=D_{12}.\end{cases}\end{equation}
For that computation we may use any representative of the isomorphism classes of
curves with absolute invariant $t$,
for example the curves $X^6+X^3+t$ or $X^5+X^3+tX$,
but of course starting from different,
non $\Q$-isomorphic curves defined over $\Q$ with the same absolute invariant
produces elliptic quotients defined in general over different fields,
that may not be isomorphic over a common field of definition.

From the expression (\ref{j_invariants_quotients}) we see that the field of moduli of the
elliptic quotients is the field $\Q(\sqrt t)$, that is the field of the rationals if
$t$ is a square or a quadratic field otherwise,
except for the single value $t=-81/700$ in the case $\Aut(C)\simeq D_8$
and the two values $t=-1/50,-4/11$ in the case $\Aut(C)\simeq D_8$,
that are non-squares but the field of moduli of the elliptic quotients is $\Q$
instead of $\Q(\sqrt t)$;
it is easily checked that these three exceptions correspond to $j$-invariants
of elliptic curves with complex multiplication.

Since the $j$-invariants in (\ref{j_invariants_quotients}) are rational or quadratic numbers
and all the orders of quadratic fields with class number one or two are known,
it is easy to find all values of $t$ for which the elliptic quotients of a curve
with absolute invariant $t$ has complex multiplication, and one gets the values
$$t=\frac{9}{100},\frac{81}{196},\frac{3969}{16900},\frac{-81}{700},\frac{1}{5},
\frac{9}{32},\frac{12}{49},\frac{81}{320},\frac{81}{325},\frac{2401}{9600},
\frac{9801}{39200},\frac{6480}{25921},\frac{194481}{777925},\frac{96059601}{384238400}$$
in the $D_8$ case and the values
$$t=\frac{4}{25},\frac{-1}{50},\frac{-4}{11},\frac{1}{20},\frac{1}{2},\frac{27}{100},
\frac{4}{17},\frac{125}{484},\frac{20}{81},\frac{256}{1025},\frac{756}{3025},\frac{62500}{250001}$$
in the $D_{12}$ case.

\paragraph{Parametrization of $\Q$-curves of degrees $2$ and $3$.}
Elkies proved that the $\Q$-curves without complex multiplication are parametrized,
up to isogeny, by the rational points of the modular curves $X^*(N)$,
quotient of the classical modular curves $X_0(N)$ by all the Atkin-Lehner involutions,
for squarefree integers $N$.
Every point of $X^*(N)(\Q)$ corresponds to a $\Q$-curve
whose isogenies to its Galois conjugates have degrees dividing $N$ up to squares.
In particular, if $N=p$ is a prime number,
the $\Q$-curves parametrized by $X^*(p)(\Q)$
are usually known as \emph{$\Q$-curves of degree $p$};
they are isogenous either to a curve defined over $\Q$
having an isogeny of degree $p$ to another curve defined over $\Q$,
or to a curve defined over a quadratic field
having an isogeny of degree $p$ to its Galois conjugate.
In the second case we say that the curve is a \emph{quadratic} $\Q$-curve.

Explicit parametrizations of $X^*(N)$ have been computed by several authors
in the case that $N$ is a prime number and $X_0(N)$ has genus zero,
and in \cite{peplario} a method is given for computing such a parametrization
that works for all curves $X^*(N)$ of genus zero or one.
For the values $N=2,3$ we follow \cite{peplario} for parametrizing $X^*(N)(\Q)$.
One may take the functions $h_0=h_0(\tau)$ or $h=h(\tau)$
$$h_0(\tau)=\begin{cases} \displaystyle\frac{\eta(\tau)^{24}}{\eta(2\tau)^{24}}+
  \displaystyle\frac{2^{12}\eta(2\tau)^{24}}{\eta(\tau)^{24}},&\ \ N=2,\\
  \phantom{\ }\\
  \displaystyle\frac{\eta(\tau)^{12}}{\eta(3\tau)^{12}}+
  \displaystyle\frac{3^6\eta(3\tau)^{12}}{\eta(\tau)^{12}},&\ \ N=3,\end{cases}\qquad\quad
  h=\begin{cases} \displaystyle\frac{h_0+104}{4(h_0-152)},&\ \ N=2\\
                  \phantom{\ }\\
                  \displaystyle\frac{h_0+42}{4(h_0-66)},&\ \ N=3.\end{cases}$$
as generators of the function field of the curve $X^*(N)$ over $\Q$,
where $\eta$ is the Dedekind $\eta$-function and $\tau$ is a variable
taking values in the upper half plane.
Then, one expresses the Galois conjugated $j$-invariants
$j(\tau)$ and $j(N\tau)$ of two elliptic curves linked by an isogeny of degree $N$
in terms of rational values of the function $h$ and obtains the following expression:
\begin{equation}\label{j_invariants_qcurves}
(j(\tau),j(N\tau))=\begin{cases}
\displaystyle\frac{2^6(3\pm10\sqrt h)^3}{(1\pm2\sqrt h)(1\mp2\sqrt h)^2},&\qquad N=2,\\
\phantom{\ }\\
\displaystyle\frac{2^8 3^3(2\pm5\sqrt h)^3(\mp\sqrt h)}{(1\pm2\sqrt h)(1\mp2\sqrt h)^3},&
\qquad N=3.\end{cases}\end{equation}
Of course, the choice of the function $h$ instead of the function $h_0$
is only motivated by aesthetic reasons,
since with this function as the parameter we obtain exactly
the same parametrization of the $j$-invariants of $\Q$-curves
than that of the $j$-invariants of elliptic quotients of curves of genus 2
expressed in terms of the absolute invariant $t$.
The identity between the expressions (\ref{j_invariants_quotients}) and
(\ref{j_invariants_qcurves}) is a direct and constructive proof of the following

\begin{theorem}
Let $E/\Qb$ be an elliptic curve with no complex multiplication.
Then $E$ is a $\Q$-curve of degree $2$ (resp. $3$)
if, and only if, it is covered by some curve of genus 2 defined over $\Q$
with group of automorphisms isomorphic to $D_8$ (resp. $D_{12}$).
\end{theorem}

\paragraph{Jacobians of $\GL_2$-type.}
We consider abelian varieties up to isogeny,
and denote by $\End(A)$ the $\Q$-algebra of endomorphisms of the abelian variety.
Let $C/\Q$ be a curve of genus 2 with $\Aut(C)=D_8$ or $D_{12}$
and let $J_C$ be its Jacobian.
It is an abelian surface defined over $\Q$ that factors
up to isogeny as the square of an elliptic curve $E$,
hence the algebra of endomorphisms $\End(J_C)$ is a 2-dimensional matrix algebra
over the field $\End(E)$, which is equal to $\Q$ or to an imaginary quadratic field
depending on whether the curve $E$ has complex multiplication or not.
Assume from now on that $E$ has no complex multiplication
(as discussed earlier, this restriction only affects the
finite number of CM $j$-invariants obtained from the given values of the parameter $t$).
Then $\End(J_C)\simeq M_2(\Q)$.
The structure of the $\Q$-algebra of endomorphisms defined over $\Q$
is given in the following

\begin{proposition}
Let $C/\Q$ be a curve of genus 2 with $\Aut(C)\simeq D_8$ or $D_{12}$
whose elliptic quotients have no complex multiplication.
Then the subalgebra $\End_\Q(J_C)$ of endomorphisms defined over $\Q$
depends only on the type of Galois action on $\Aut(C)$
and is reduced to $\Q$ except in the cases listed in the following tables
$$\begin{array}{|c||c|c|c|c|}
\hline
\text{Type} & \Q-basis & \End_\Q(J_C) & \Q(j) & \text{condition} \\
\hline\hline
I & 1,U,V,UV & M_2(\Q) & \Q & \\
\hline
C_2^A & 1,V & \Q(\sqrt{-1}) & \Q & \\
\hline
C_2^B & 1,U & \Q\times\Q & \Q & \\
\hline
C_2^C & 1,U(1+V) & \Q(\sqrt2) & \Q(\sqrt t) & (t,2)=1 \\
\hline
C_4 & 1,V & \Q(\sqrt{-1}) & \Q(\sqrt t) & (w,-1)=1 \\
\hline
\end{array}$$
\begin{center}{\bf Table 3.} Nontrivial endomorphisms of $J_C$ defined over $\Q$.
  Case $\Aut(C)\simeq D_8$.\end{center}

$$\begin{array}{|c||c|c|c|c|}
\hline
\text{Type} & \Q-basis & \End_\Q(J_C) & \Q(j) & \text{condition} \\
\hline\hline
I & 1,U,V,UV & M_2(\Q) & \Q & \\
\hline
C_2^A & 1,V & \Q(\sqrt{-3}) & \Q(\sqrt t) & (t,-3)=1 \\
\hline
C_2^B & 1,U & \Q\times\Q & \Q & \\
\hline
C_2^C & 1,U(1+2V) & \Q(\sqrt3) & \Q(\sqrt t) & (t,3)=1 \\
\hline
C_3 & 1,V & \Q(\sqrt{-3}) & \Q & \\
\hline
C_6 & 1,V & \Q(\sqrt{-3}) & \Q(\sqrt t) & (t,-3)=1 \\
\hline
\end{array}$$
\begin{center}{\bf Table 4.} Nontrivial endomorphisms of $J_C$ defined over $\Q$.
  Case $\Aut(C)\simeq D_{12}$.\end{center}
The column $\Q(j)$ on the tables gives the field of definition of the $j$-invariant
of the quotient elliptic curves, that can be either $\Q$ or a quadratic field,
and the last column gives the obstruction in $\Br_2(\Q)$ to the existence of a curve $C$
with the given absolute invariant and Galois action on automorphisms.
\end{proposition}

\begin{proof}
There is a natural embedding of $\Aut(C)$ in $\Aut(J_C)\subset\End(J_C)$
that respects the field of definition of the morphisms.
Since the automorphisms of the jacobian corresponding to $1,U,V,UV\in\Aut(C)$
are independent over $\Q$, they are a basis of $\End(J_C)$ as a $\Q$-vector space,
and we can find the Galois action on arbitrary endomorphisms from
the knowledge of the action on that basis.
Then the tables are built just by checking case by case.
The last column is the obstruction of corollaries
\ref{theorem_obstruction_D8} and \ref{theorem_obstruction_D12}
expressed in terms of the parameter $t$ for the $D_{12}$ case
and in terms of the parameters $t$ and $w$ for the $D_8$ case.
\end{proof}

\paragraph{Quotient $\Q$-curves completely defined over certain fields.}
Let $E/K$ be a $\Q$-curve with no complex multiplication that is
completely defined over a Galois number field $K$.
For every $\sigma\in\Gal(K/\Q)$ choose an isogeny $\phi_\sigma:{}^\sigma E\to E$.
Then, for every pair of elements $\sigma,\tau\in\Gal(K/\Q)$, the map
$\phi_\sigma\circ{}^\sigma\phi_\tau\circ\phi_{\sigma\tau}^{-1}:E\to E$
is an isogeny and it can be identified with a nonzero rational number.
We may define cocycles of $\Gal(K/\Q)$ with values in $\Q^*$ and in $\{\pm1\}$ by
$$c(\sigma,\tau)=\phi_\sigma\circ{}^\sigma\phi_\tau\circ\phi_{\sigma\tau}^{-1},\qquad
  c^\pm=\sign(c).$$
The cohomology classes of these cocycles,
which we denote in the same way,
are invariants of the $K$-isogeny class of the curve.
The decomposition $\Q^*=\{\pm1\}\times P$,
with $P$ the multiplicative group of positive numbers,
induces a decomposition $H^2(K/\Q,\Q^*)\simeq H^2(K/\Q,\{\pm1\})\times H^2(K/\Q,P)$;
the first component of $c$ is just $c^\pm$ and the second one depends only on
the isogeny class of the curve (it is in fact an invariant of the $\Qb$-isogeny class).

\begin{lemma}\label{lemma_Q-twist}
Let $E$ and $E'$ be isomorphic elliptic curves both completely defined over the
Galois number field $K/\Q$.
If the two cohomology classes $c^\pm,{c'}^\pm\in H^2(K/\Q,\{\pm1\})$ are the same
then the two curves are isomorphic over a quadratic extension $K(\sqrt d)$
for a rational number $d$.
\end{lemma}

\begin{proof}
For every elliptic curve $E$ defined over a field $k$,
and element $\gamma\in k^*$,
we denote by $E_\gamma$ the twisted curve defined over $k$
and isomorphic to $E$ over $k(\sqrt\gamma)$,
which is completely determined by $\gamma$ modulo squares.
If $Y^2=f(X)$ is a Weierstrass equation for $E$ then $\gamma Y^2=f(X)$
is an equation for $E_\gamma$.

The two given curves, being isomorphic and both defined over $K$,
must be a twist of each other by some element $\gamma\in K^*$, say $E'=E_\gamma$.

For every isogeny $\phi_\sigma:{}^\sigma E\to E$ let $\lambda_\sigma\in K^*$
be the number defined by the identity
$\phi_\sigma^*(\omega)=\lambda_\sigma{}^\sigma\omega$,
for an invariant differential $\omega$ of $E/K$.
Then $c(\sigma,\tau)=\lambda_\sigma{}^\sigma\lambda_\tau\lambda_{\sigma\tau}^{-1}$.
The same construction for the curve $E'$ produces elements $\mu_\sigma\in K^*$
with $c'(\sigma,\tau)=\mu_\sigma{}^\sigma\mu_\tau\mu_{\sigma\tau}^{-1}$,
and the fact that $E'=E_\gamma$ implies the identity
$\mu_\sigma^2={}^\sigma\gamma\gamma^{-1}\lambda_\sigma^2$.
The hypotheses that the sign components $c^\pm$ and ${c'}^\pm$ are the same
and that the curves are isomorphic (hence isogenous)
imply that the two cohomology classes $c,c'$ are the same.
Then, if we define $b_\sigma=\mu_\sigma\lambda_\sigma^{-1}\in K^*$,
the map $(\sigma,\tau)\mapsto b_\sigma{}^\sigma b_\tau b_{\sigma\tau}^{-1}$
is a 2-cocycle with values in $\Q^*$ and trivial cohomology class;
changing our choice of isogenies $\phi_\sigma$ by rational numbers
accordingly to a boundary map for this 2-cocycle,
we may assume that it is the trivial map.
Then, the map $\sigma\mapsto b_\sigma$ is a one-cocycle of $\Gal(K/\Q)$ with values in $K^*$
and, by Hilbert's Theorem 90, there exist an element $\alpha\in K^*$
such that $b_\sigma={}^\sigma\alpha\cdot\alpha^{-1}$.
The element $d=\gamma/\alpha^2$ is then fixed by every $\sigma\in\Gal(K/\Q)$.
Hence it is a rational number and the two given curves are isomorphic over $K(\sqrt d)$.
\end{proof}

Let $B=\Res_{K/\Q}(E)$ be the abelian variety of dimension $[K:\Q]$ defined over $\Q$
obtained by restriction of scalars.
It is of interest to characterize the situation where the
variety $B$ factors, up to $\Q$-isogeny, as a product of abelian varieties of $\GL_2$-type,
since, assuming the modularity of the $\Q$-curve $E$,
this is equivalent to the fact that the $L$-series of the curve $E$ over the field $K$
is a product of $L$-series of modular forms for congruence subgroups $\Gamma_1(N)$.
In \cite[Proposition 5.2]{quer} it is shown that the necessary and sufficient condition
for that fact is that the extension $K/\Q$ is abelian
and the cohomology class $c^\pm\in H^2(K/\Q,\{\pm1\})$
belongs to the subgroup $\Ext(K/\Q,\{\pm1\})$,
consisting of classes represented by symmetric 2-cocycles.

For the curves of genus 2 that we study in this paper
there are several cases of elliptic quotients that have the property discussed in the
previous paragraph.
The following theorem is a strong converse of that fact
since it gives sufficient conditions for $\Q$-curves completely defined over a field $K$
to be covered by a curve of genus 2 defined over $\Q$
with covering defined over $K$.

\begin{theorem}\label{theorem_quotients_completely_defined}
\begin{enumerate}
\item Let $E/K$ be a quadratic $\Q$-curve of degree 3
completely defined over a quadratic field $K$.
Then there is a curve $C/\Q$ of genus 2 with $\Aut(C)\simeq D_{12}$
and a covering $C\to E$ defined over $K$.

\item Let $E/K$ be a quadratic $\Q$-curve of degree 2
completely defined over a quadratic field $K$.
Then there is a curve $C/\Q$ of genus 2 with $\Aut(C)\simeq D_8$
and a covering $C\to E$ defined over $K$ if, and only if,
the class of $c^\pm\in H^2(K/\Q,\{\pm1\})$ is trivial.

\item Let $E/K$ be a quadratic $\Q$-curve of degree 2
completely defined over a cyclic quartic field $K$
with $j$-invariant in the quadratic subfield.
Then there is a curve $C/\Q$ of genus 2 with $\Aut(C)\simeq D_8$
and a covering $C\to E$ defined over $K$.
\end{enumerate}
\end{theorem}

\begin{proof}
Let $E/K$ be a $\Q$-curve as in 1.
Let $j(E)$ be given by the parametrization (\ref{j_invariants_qcurves})
from some rational value $h\in\Q$. Then $K=\Q(\sqrt h)$.

The group $H^2(K/\Q,\{\pm1\})$ is isomorphic to $\Z/2\Z$;
we denote by $c_\varepsilon$ the nontrivial element.
In the interpretation of that group as classifying group extensions of $\Gal(K/\Q)$
with kernel the group $\{\pm1\}$, the trivial element corresponds to the extension
realized by the group $V_4$ and the nontrivial one $c_\varepsilon$
to that by the group $C_4$.

By the formula of \cite[Theorem 3.1]{quer}
the inflation $\Inf(c^\pm)\in H^2(G_\Q,\{\pm1\})\simeq\Br_2(\Q)$
corresponds to the element $(h,3)$.
By the general theory of embedding problems in Galois theory,
the inflation $\Inf(c_\varepsilon)$ is the obstruction to embed the quadratic field $K$
into a cyclic extension of degree 4, which is the element $(h,-1)\in\Br_2(\Q)$.

Assume that $c^\pm$ is trivial, then $(h,3)=1$ and by theorem
(\ref{theorem_obstruction_D12})
there exists a curve $C/\Q$ of genus 2 with $\Aut(C)\simeq D_{12}$
having absolute invariant $t=h$,
and with Galois action on $\Aut(C)$ of type $C_2^C$ defined over $K=\Q(\sqrt t)$.
Let $E'/K$ be a quotient of $C$ by an hyperelliptic involution
with the same $j$-invariant than $E$,
which exists due to the identity between (\ref{j_invariants_qcurves})
and (\ref{j_invariants_quotients}).
The map $C\to E'$ is defined over the field $K$ and by the universal property of
restriction of scalars there is an isogeny $J_C\to\Res_{K/\Q}(E')$ defined over $\Q$.
By the computation of the endomorphisms defined over $\Q$ given in Table 3
we know that $\End_\Q(J_C)\simeq\Q(\sqrt3)$.
It follows that $\End_\Q(\Res_{K/\Q}(E'))\simeq\Q(\sqrt3)$ and,
by \cite[Theorem 5.4]{quer}, the 2-cocycle class ${c'}^\pm$ corresponding to the
curve $E'$ must be the trivial class in $H^2(K/\Q,\{\pm1\})$.

The curves $E$ and $E'$ are isomorphic, both completely defined over $K$
and the cohomology classes $c^\pm$ and ${c'}^\pm$ are the same (the trivial class).
By lemma \ref{lemma_Q-twist} they are a twist of each other
over some field $K(\sqrt d)$ for an element $d\in\Q^*$.
Let $C_d$ be the hyperelliptic twist of the curve $C$ corresponding to this element.
Then the curve $E=E'_d$ is one of the quotients of $C_d$ by an involution
and the map $C_d\to E$ is defined over $K$.

Assume now that $c^\pm$ is the nontrivial element $c_\varepsilon$.
Then their inflations coincide, $(h,3)=(h,-1)$ and this implies $(h,-3)=1$.
By theorem \ref{theorem_obstruction_D12} there exists a curve $C/\Q$
of genus 2 with $\Aut(C)\simeq D_{12}$ having absolute invariant $t=h$,
and with Galois action on $\Aut(C)$ of type $C_2^A$ defined over $K=\Q(\sqrt t)$.
Now, the argument follows in exactly the same way than in the previous case,
with the construction of an elliptic curve $E'/K$ that is a quotient of $C$,
with $C\to E'$ defined over $K$,
hence $\End_\Q(C)\simeq\End_\Q(\Res_{K/\Q}(E'))\simeq\Q(\sqrt{-3})$ by Table 3,
with cocycle ${c'}^\pm$ equal to $c_\varepsilon$ by \cite[Theorem 5.4]{quer},
hence a twist of $E$ over a field $K(\sqrt d)$ with $d\in\Q^*$,
from which we find a curve $C_d$ with a map $C_d\to E$ defined over $K$.

The part 2 can be proved in the same way than part 1 in the case of trivial $c^\pm$,
just by remarking that under the condition $(h,2)=1$ there is a curve of genus 2
with absolute invariant $t=h$
and with Galois action on $\Aut(C)$ of type $C_2^C$ defined over $K=\Q(\sqrt t)$.
The only if part of the statement comes from the fact that when $c^\pm=c_\varepsilon$,
then by \cite[Theorem 5.4]{quer} $\End_\Q(\Res_{K/\Q}(E))\simeq\Q(\sqrt{-2})$
and this field does never happen to be the field of endomorphisms of
the jacobian of a curve of genus 2 in the case $D_8$, as can be checked in Table 3.

Now we consider the third part of the theorem. Let $E/K$ be such a curve.
Since $j(E)$ belongs to a quadratic field,
it can be obtained by the parametrization (\ref{j_invariants_qcurves})
from a rational value $h\in\Q$
with $\Q(\sqrt h)$ being the quadratic subfield of $K$.
We may assume that $K$ is constructed as in lemma \ref{lemma_expression_field_D8}
with $u=h$ and some rational numbers $w\neq0$ and $z$.

The group $H^2(K/\Q,\{\pm1\})$ is again isomorphic to $\Z/2\Z$ and
we denote by $c_\varepsilon$ the nontrivial element.
In terms of group extensions the two elements correspond to extensions by
the groups $C_2\times C_4$ and $C_8$.
The inflation of $c^\pm$ to $G_\Q$ is $(h,2)$.
It is an exercise in the theory of Galois embedding problems to compute
the obstruction to embed the given cyclic quartic field $K$
into a cyclic extension of degree 8,
and one obtains the element $(h,2)(w,-1)\in\Br_2(\Q)$.

If $c^\pm$ is trivial, then $(h,2)=1$ and we apply the previous case
obtaining a curve $C$ and a map $C\to E$ defined over the quadratic subfield $\Q(\sqrt h)$,
and a fortiori over $K$.
The case $c^\pm=c_\varepsilon$ is more interesting.
Then by inflation we obtain the identity $(h,2)=(h,2)(w,-1)$ that implies $(w,-1)=1$,
and by theorem \ref{theorem_obstruction_D8}
there exists a curve $C/\Q$ of genus 2 with $\Aut(C)\simeq D_8$
having absolute invariant $t=h$,
and with Galois action on $\Aut(C)$ of type $C_4$ defined over $K$.
Let $E'/K$ be a quotient of $C$ by an hyperelliptic involution
with the same $j$ invariant that $E$,
which exists due to the identity between (\ref{j_invariants_qcurves})
and (\ref{j_invariants_quotients}).
The map $C\to E'$ is defined over the field $K$ and by the universal property of
restriction of scalars there is a nontrivial morphism $J_C\to\Res_{K/\Q}(E')$ defined over $\Q$.
By the computation of the endomorphisms defined over $\Q$ in Table 3
we know that $\End_\Q(J_C)\simeq\Q(\sqrt{-1})$;
it follows that $\Res_{K/\Q}(E')$ has a 2-dimensional factor up to isogeny
of $\GL_2$-type with $\Q$-endomorphisms $\Q(\sqrt{-1})$ and,
by \cite[Theorem 5.4]{quer}, the 2-cocycle class ${c'}^\pm$ corresponding to the
curve $E'$ must be the class $c_\varepsilon$ in $H^2(K/\Q,\{\pm1\})$.

Now we can use the same argument than in the proof of the first part:
the curves $E$ and $E'$
are twists of each other over a quadratic extension $K(\sqrt d)$ for some $d\in\Q^*$
the hyperelliptic twist $C_d$ has the curve $E$ as one of its quotients
by a non-hyperelliptic involution and the map $C_d\to E$ is defined over $K$.
\end{proof}

The results of this Theorem interpret the fields that appear in Tables 3 an 4
as fields endomorphisms of jacobians of curves of genus 2,
but there are still some remarks to be done that improve the understanding
of the maps $C\to E$ covering elliptic $\Q$-curves:

For the elliptic curves $E/\Q$ of the third part of the theorem
the abelian variety $\Res_{K/\Q}(E)$ factors over $\Q$ up to isogeny
as a product of two abelian surfaces of $\GL_2$-type
that are not isogenous over $\Q$ but over the quadratic field $\Q(\sqrt t)$
contained in the quartic field $K$.
If $C\to E$ is a covering defined over $K$ then the hyperelliptic twist $C_t$
is non-isomorphic to $C$ over $\Q$
(cf. proposition \ref{proposition_hyperelliptic_twists})
and it also covers the same curve $E$ over the field $K$,
since the twist $E_t$ is $K$-isomorphic to $E$ because $\Q(\sqrt t)\subset K$.
The two abelian surfaces in which $\Res_{K/\Q}(E)$ factors
are up to $\Q$-isogeny the two jacobians $J_C$ and $J_{C_t}$.

If $(t,-2)=1$ then there are quadratic $\Q$-curves $E$ of degree 2
completely defined over the quadratic field $K=\Q(\sqrt t)$,
and the abelian surfaces $B=\Res_{K/\Q}(E)$ are of $\GL_2$-type
with $\End_\Q(B)\simeq\Q(\sqrt{-2})$,
but these curves never appear as covered by genus 2 curves $C/\Q$ with $\End(C)\simeq D_8$,
with $C\to E$ defined over $K$,
and these abelian varieties are never isogenous over $\Q$ to the jacobian of one
of the curves of genus 2 under study.
This is the reason of the missing of the field $K(\sqrt{-2})$ in Table 3.
Indeed, from the identity $(t,-2)=1$ and applying theorem \ref{theorem_obstruction_D8},
there exist curves $C/\Q$ of genus 2 with absolute invariant $t$
and group of automorphisms $D_8$
defined over the biquadratic field $K=\Q(\sqrt t,\sqrt{-2})$
(and also over many other biquadratic fields).
The quotient elliptic curves $E$ and the coverings $C\to E$
are defined over that biquadratic field.
Even if we compose with an isomorphism $E\to E'$,
with $E'$ completely defined over $\Q(\sqrt t)$, which always exists,
the covering $C\to E\to E'$ can never be defined over the quadratic field.
Moreover, and even worse, the four dimensional abelian variety $\Res_{K/\Q}(E)$
has no factors of $\GL_2$-type since the cocycle $c^\pm$
does not belong to $\Ext(K/\Q,\{\pm1\})$,
and this prevents the jacobian of the curve $C$ from having that property.


\vfill
\noindent
\begin{tabular}{ll}
\begin{minipage}[t]{8cm}
Gabriel Cardona\\
Dept. Ci\`encies Matem\`atiques i Inf.\\
Universitat de les Illes Balears\\
Ed. Anselm Turmeda, Campus UIB\\
Carretera Valldemossa, km. 7.5\\
07071 -- Palma de Mallorca, Spain\\
\verb·gabriel.cardona@uib.es·
\end{minipage}
&
\begin{minipage}[t]{8cm}
Jordi Quer\\
Dept. Matem\`atica Aplicada II\\
Universitat Polit\`ecnica de Catalunya\\
Ed. U, Campus Sud\\
Pau Gargallo, 5\\
08028 -- Barcelona, Spain\\
\verb·quer@ma2.upc.es·
\end{minipage}
\end{tabular}
\end{document}